\documentclass[11pt,a4paper,reqno,twoside]{article}
\usepackage{amsmath,amsthm,amsfonts,amssymb}
\usepackage{mathrsfs}
\usepackage{enumitem}
\usepackage[utf8]{inputenc}
\usepackage[T1]{fontenc}
\usepackage{tgtermes}
\usepackage{hyperref}
\hypersetup{  
   bookmarks=true,
   backref=true,
   pagebackref=false,
   colorlinks=true,
   linkcolor=blue,
   citecolor=red,
   urlcolor=blue
}
\setlist{nolistsep}

\usepackage{subfig}
\usepackage[english]{babel}
\usepackage{amsbsy,amscd,fancyhdr,graphicx,psfrag,fancybox,indentfirst,color}
\usepackage{graphics,epsfig}
\usepackage{mysty}
\usepackage{fullpage}
\usepackage{graphicx}

\newcommand*\samethanks[1][\value{footnote}]{\footnotemark[#1]}

\title{Continuum limit of $p$-Laplacian evolution problems on graphs: $L^q$ graphons and sparse graphs}         
\date{}
\author{Imad El Bouchairi\thanks{Normandie Univ, ENSICAEN, UNICAEN, CNRS, GREYC, France.} \and Jalal M. Fadili\samethanks \and Abderrahim Elmoataz\samethanks[1]
}

\begin{document}

\maketitle

\begin{abstract}
In this paper we study continuum limits of the discretized $p$-Laplacian evolution problem on sparse graphs with homogeneous Neumann boundary conditions. This extends the results of \cite{Hafienne_FE_Nonlocal_plaplacian_evolu} to a far more general class of kernels, possibly singular, and graph sequences whose limit are the so-called $L^q$-graphons. More precisely, we derive a bound on the distance between two continuous-in-time trajectories defined by two different evolution systems (i.e. with different kernels, second member and initial data). Similarly, we provide a bound in the case that one of the trajectories is discrete-in-time and the other is continuous. In turn, these results lead us to establish error estimates of the full discretization of the $p$-Laplacian problem on sparse random graphs. In particular, we provide rate of convergence of solutions for the discrete models to the solution of the continuous problem as the number of vertices grows.
\end{abstract}

\begin{keywords}
Nonlocal diffusion; $p$-Laplacian; sparse graphs; $L^q$-graphons.
\end{keywords}

\begin{AMS}
65N12, 34A12, 45G10, 05C90.
\end{AMS}


\section{Introduction}
\subsection{Problem formulation}

Our main goal in this paper is to study discretization of the following nonlinear diffusion problem on graphs, which we call the nonlocal $p$-Laplacian problem with homogeneous Neumann boundary conditions: 
\begin{equation}\tag{\textrm{$\mathcal{P}$}}
\begin{cases}
\frac{\partial }{\partial t} u(\bx,t) = \int_{\O} \K(\bx,\by) \abs{u(\by,t) - u(\bx,t) }^{p-2} (u(\by,t) - u(\bx,t) ) d\by + f(\bx,t), & \bx \in \O, t > 0,
\\
u(\bx,0)= g(\bx), \quad \bx \in \O,
\end{cases}
\label{neumann}
\end{equation}
where $p \in [1,+\infty[$, $\O \subset \R^d$ is a bounded domain, $d \geq 1$, without loss of generality $\O = [0,1]^d$, and $\K: \R^d \times \R^d \to \R$ is the kernel function. In particular, in the setting of graphs, $d=1$ and it will be seen that $\K$ is the limit object for some convergent graph sequence $\{G_n\}_n, n \in \N$, whose meaning and form will be specified in the sequel. Throughout, we assume that
\begin{enumerate}[label=({\textbf{H.\arabic*}})]
\item $\K$ is a nonnegative measurable function. \label{assum:Kpos}
\item $\K$ is symmetric, i.e., $\K(\bx,\by)=\K(\by,\bx)$. \label{assum:Ksym}
\item $\sup_{\bx \in \O} \int_{\O}\K(\bx,\by) d\by < +\infty$ . \label{assum:KL1}
\end{enumerate}
By \ref{assum:Ksym}, it is straightforward to see that 
\[
\sup_{\bx \in \O} \int_{\O}\K(\bx,\by) d\by = \sup_{\by \in \O} \int_{\O}\K(\bx,\by) d\bx ,
\]
and thus, \ref{assum:KL1} is equivalent to
\[
\sup_{\by \in \O} \int_{\O} \K(\bx,\by) d\bx < +\infty .
\]
When the kernel is such that $\K(\bx,\by) = \J(\bx-\by)$, where $\J: \R^d \to \R$, then \ref{assum:Kpos}, \ref{assum:Ksym} and \ref{assum:KL1} read:
\begin{enumerate}[label=({\textbf{H'.\arabic*}})]
\item $\J$ is nonnegative and measurable. \label{assum:Jpos}
\item $\J$ is symmetric, i.e., $\J(-\bx)=\J(\bx)$. \label{assum:Jsym}
\item $\int_{\O-\O}\J(\bx) d\bx < +\infty$ . \label{assum:JL1}
\end{enumerate}
Recall that $\O-\O$ is the Minkowski sum of $\O$ and $-\O$. In the case $\O = [0,1]^d$, we obviously have $\O-\O=[-1,1]^d$.

The motivation behind this work is that partial differential equations (PDEs) involving the nonlocal $p$-Laplacian operator have become more and more popular both in the setting of Euclidean domains and on discrete graphs, as the $p$-Laplacian problem has been possessing many important features shared by many practical problems in mathematics, physics, engineering, biology, and economy, such as continuum mechanics, phase transition phenomena, population dynamics, see~\cite{mazon, phase, waves, effect, like, patterns, interfacial} and references therein. Some closely related applications can be found in image processing, computer vision and machine learning~\cite{cluster, image, regularization, elliptic}.

The main goal of this paper is to revisit and extend the work of \cite{Hafienne_FE_Nonlocal_plaplacian_evolu} by removing important limiting assumptions made there on the kernel $\K$ and the initial condition $g$. In turn, this will allow us to establish consistency estimates of the fully discretized $p$-Laplacian problem for singular kernels or on sparse graphs whose limits are known not to be bounded graphons (see Section~\ref{graphlimits} for details).


\subsection{Contributions}
In this work, we propose a far-reaching generalization of the results in \cite{Hafienne_FE_Nonlocal_plaplacian_evolu} to a much more general class of kernels and initial data. In particular, we are able to consider unbounded initial data, the case $p=1$, and most importantly singular kernels, which in turn will allow to handle sparse graph sequences whose limit are the so-called $L^q$-graphons \cite{Borgs18partII,Borgs19partI}. On such graphs sequences, we will quantitatively analyse evolution problems and their continuum limit. We will also consider the case $p=1$ which was not handled in \cite{Hafienne_FE_Nonlocal_plaplacian_evolu}.

More precisely, for $p \in [1,+\infty[$, we first provide a general error bound between the solutions of two evolution problems of the form~\eqref{neumann} governed by two kernels and two initial data. This will then be applied to study consistency of numerical solutions regarded as a discretization of~\eqref{neumann} with space discretization of the kernel $J$ and initial data $g$, and forward and backward Euler discretization in time. We eventually apply these results to evolution problems defined on graph sequences whose limit are the so-called $L^q$-graphons. 

\subsection{Relation to prior work}
The kernels and initial data considered here are beyond reach of the approach developed in \cite{Hafienne_FE_Nonlocal_plaplacian_evolu}, and have not been considered in the literature to the best of our knowledge. Moreover, our error bounds are directly stated in $L^2(\O)$ and not in $L^p(\O)$ as done in this previous work. Our proof is also simpler, more elegant and the argument is made more transparent. This argument will allow us to handle the case $p=1$. More importantly, some limiting assumptions on the kernel and the initial data made in \cite{Hafienne_FE_Nonlocal_plaplacian_evolu} are removed and replaced by much less stringent ones. This allows in particular to cover a far larger class of kernels (including singular ones), and also sparse graph sequences that were not handled in that previous work.

Another related work is that in~\cite{Medvedevsparse,medvedevkuramotosparse}. In these papers, the authors focused on a nonlinear heat equation on sparse graphs, where Lipschitz-continuity of the operator is of paramount importance. This assumption was essential to prove well-posedness (existence and uniqueness follow immediately from the contraction principle), as well as to study the consistency in $L^2(\O)$ of the spatial semi-discrete approximation. The nonlocal $p$-Laplacian evolution problem considered here is much more general and cannot be covered by the approach of those previous papers because the lack of Lipschitzianity raises several challenges (including for well-posedness and error estimates). Unlike those previous works, we also consider both the semi-discrete and fully-discrete versions with both forward and backward Euler approximations, that we fully characterize, and develop novel proof techniques.

\subsection{Paper organization}
The rest of the paper is organized as follows. In Section \ref{prelim} we start by reviewing some basic notations and recall some preliminary material necessary to our exposition. In Section~\ref{graphlimits} , we provide some prerequisites on $L^q$-graphons and sparse $\K$-random graph models that we are going to deal with. Section~\ref{existenceuniqueness} is devoted to study the well-posedness of the problem~\eqref{neumann}. In Section~\ref{sec:consistencecontinuous}, we study stability of the problem~\eqref{neumann} with respect to sequences of kernels $\K$, initial data $g$ and second member $f$. Error bounds for the semi-discete (i.e., space discretization of $(\K,g,f)$) problem are established in Section~\ref{semidiscrete}, and those  for the fully discrete (time and space discretization) problem with forward and backward Euler time-discretization are provided in Section~\ref{totallydiscrete}. Section~\ref{sec:applicationgraphs} is devoted to applying these results to fully discretized problems on sparse random graph models.

\section{Preliminaries}
\label{prelim}

\subsection{Basic notations}
For $q \in [1,+\infty]$, and $S \subset \R^d$, $L^q(S)$ is the standard Banach space of Lebesgue $q$-integrable functions on $S$. For a function $F: S \times S \to \R$, we define the $L^{\infty,q}(S^2)$-norm as 
\[
\norm{F}_{L^{\infty,q}(S^2)} \eqdef \sup_{\bx \in S} \norm{F(\bx,\cdot)}_{L^q(S)} .
\]
If $F$ is symmetric, then
\[
\norm{F}_{L^{\infty,q}(S^2)} = \sup_{\by \in S} \norm{F(\cdot,\by)}_{L^q(S)} .
\] 
$L^{\infty,q}(S^2)$ is the space of functions on $S^2$ of bounded $L^{\infty,q}(S^2)$)-norm, which is of course a Banach space.

Throughout the paper, we will often use Fubini's theorem without explicitly referring to it.

$C([0,T];L^q(S))$ denotes the space of functions $u: S \times [0,T] \to \R$ which are uniformly continuous in time with values $u(t,\cdot)$ in $L^p(S)$. $C([0,T];L^q(S))$ is naturally endowed with the norm $\norm{u}_{C([0,T];L^q(S))} \eqdef \sup_{t \in [0,T]} \norm{u(\cdot,t)}_{L^q(S)}$. Moreover, $L^1([0,T];L^q(S))$ is the space of functions $u: S \times [0,T] \to \R$ such that $\norm{u}_{L^1([0,T];L^q(S))} \eqdef \int_{0}^T \norm{u(\cdot,t)}_{L^q(S)} dt < +\infty$.

We define the $q$-norm on $\R^{m}$, $q \in [1, +\infty[$, as
\[
\norm{\mb v}_{q}= \pa{\frac{1}{m} \sum_{i \in [m]} \aabs{\mb v_{i}}^q}^{\frac{1}{q}} ,
\]
with the usual adaptation for $q=+\infty$. $(\cdot)_+$ is the positive part function on $\R$. 

For a set-valued operator $A: X \to 2^X$ on a Banach space $X$, its domain and range are respectively 
\[
\dom(A) = \ens{u \in X}{A u \neq \emptyset} \qandq \range(A) = A(X) .
\]

\subsection{Projector and injector}
Let $n \in \N^*$ and denote the multi-index ${\bi} = (\bi_1,\bi_2,\ldots,\bi_d) \in [n]^d$. Partition $\O$ into cells (hypercubes)
\[
\Oni \eqdef \ens{ \prod_{k=1}^d ]\bx_{\bi_k-1}, \bx_{\bi_k}]}{{\bi} \in [n]^d}
\] 
of size $h_{\bi} \eqdef \mOi$, and maximal mesh size 
\[
\deltan \eqdef \max\limits_{{\bi} \in [n]^d} \max\limits_{k \in [d]}(\aabs{\bx_{\bi_k} - \bx_{\bi_k-1}}) .
\] 
When the cells are equispaced, then $h_{\bi}=1/n^{d}$. 

We consider the operator $\projn: L^1(\O) \to \R^{n^d}$
\begin{equation}\label{eq:projn}
(\projn u)_{\bi} \eqdef \frac{1}{h_{\bi}} \int_{\On_{\bi}} u(\bx) d\bx.
\end{equation}
This operator can be also seen as a piecewise constant projector of $u$ on the space of discrete functions. For simplicity, and with a slight abuse of notation, we keep the same notation for the projector $\projn: L^1(\O^2) \to \R^{n^d \times n^d}$.

Our aim is to study the relationship between solutions of discrete approximations and the solution of the continuum model. It is then convenient to introduce an intermediate model which is the continuum extension of the discrete solution. Towards this goal, we consider the piecewise constant injector $\injn$ of a vector $\mb v \in \R^{n^{d}}$ into $L^2(\O)$ defined as
\begin{equation}\label{eq:injn}
\injn \mb v (\bx) \eqdef \sum_{{\bi} \in [n]^d} \mb v_{{\bi}} \indOi(\bx),
\end{equation}
where we recall that $\chi_{\mathcal{C}}$ is the characteristic function of the set $\mathcal{C}$, i.e., takes $0$ on $\mathcal{C}$ and $1$ otherwise. 

It is immediate to see that the operator $\injn\projn$ is the orthogonal projector on the subspace \linebreak$\Span \ens{\indOi}{{\bi} \in [n]^d}$ of $L^1(\O)$. In turn, $\injn\projn u$ is the the piecewise constant approximation of $u$.

\begin{lem}
\label{lem:projinj}
For a function $u \in L^q(\O)$, $q \in [1, + \infty]$, we have 
\begin{equation}
\norm{\injn\projn u}_{L^q(\O)} \leq \norm{u}_{L^q(\O)} .
\label{proppvn}
\end{equation}
For a function $\K \in L^{\infty,q}(\O^2)$, $q \in [1, + \infty]$, we have
\begin{equation}
\norm{\injn\projn \K}_{L^{\infty,q}(\O^2)} \leq \norm{\K}_{L^{\infty,q}(\O)} .
\label{proppvnbi}
\end{equation}
\end{lem}

\bpf{}
We prove \eqref{proppvnbi} as \eqref{proppvn} is a consequence of it. Let $\bK = \projn K$.
We have, $\forall \bx \in \O$,
\begin{align*}
\int_{\O} \injn\projn |\K(\bx,\by)|^q d\bx d\by 
&=\int_{\O} \sum_{\bi,\bj} |\bK_{\bi,\bj}|^q \indOi(\bx)\indOj(\by) d\by \\
&=\sum_{\bi} \pa{\sum_{\bj} \int_{\Onj} |\bK_{\bi,\bj}|^q d\by} \indOi(\bx) \\
&=\sum_{\bi} \pa{\sum_{\bj} h_\bj \aabs{\frac{1}{h_\bi h_\bj}\int_{\Oni \times \Onj} \K(\bx',\by') d\bx'd\by'}^q} \indOi(\bx) \\
&\leq\sum_{\bi} \pa{\sum_{\bj} \frac{1}{h_\bi}\int_{\Oni \times \Onj}|\K(\bx',\by')|^q d\bx'd\by'} \indOi(\bx) \\
&=\sum_{\bi} \pa{\frac{1}{h_\bi}\int_{\Oni} \pa{\sum_{\bj} \int_{\Onj}|\K(\bx',\by')|^q d\by'} d\bx'}\indOi(\bx) \\
&=\sum_{\bi} \pa{\frac{1}{h_\bi}\int_{\Oni} \pa{\int_{\O}|\K(\bx',\by')|^q d\by'} d\bx'}\indOi(\bx) \\
&\leq \norm{\K}_{L^{\infty,q}(\O^2)}^q \sum_{\bi} \indOi(\bx) = \norm{\K}_{L^{\infty,q}(\O^2)}^q .
\end{align*}
Taking the supremum on the left-hand side yields the bound.
\epf{}

\subsection{Lipschitz spaces}
\label{subsec: Lipspaces}
Fro $N \in \N^*$, let $S$ be a compact subset of $\R^N$. We introduce the Lipschitz spaces $\Lip(s,L^q(S))$, $q \in [1,+\infty]$, which contain functions with, roughly speaking, $s$ "derivatives" in $L^q(S)$~\cite[Ch.~2, Section~9]{devorelorentz93}. 
\begin{defi}\label{def:lipspaces}
For $F \in L^q(S)$, $q \in [1,+\infty]$, we define the (first-order) $L^q(S)$ modulus of smoothness by
\begin{equation}
\omega(F,h)_q \eqdef
 \sup_{\bz \in \R^d, |\bz| < h} \pa{\int_{\bx, \bx + \bz \in S}\aabs{F(\bx + \bz)-F(\bx)}^q d \bx}^{1/q} .
\label{modsmooth}
\end{equation}
The Lipschitz spaces $\Lip(s,L^q(S))$ consist of all functions $F$ for which
\[
\aabs{F}_{\Lip(s,L^q(S))} \eqdef\sup_{h > 0} h^{-s} \omega(F,h)_q < +\infty .
\]
\end{defi}
We restrict ourselves to values $s \in ]0,1]$ since for $s > 1$, only constant functions are in $\Lip(s,L^q(S))$. It is easy to see that $\aabs{F}_{\Lip(s,L^q(S))}$ is a semi-norm. $\Lip(s,L^q(S))$ is endowed with the norm
\[
\norm{F}_{\Lip(s,L^q(S))} \eqdef \norm{F}_{L^q(S)} + \aabs{F}_{\Lip(s,L^q(S))} .
\]
The space $\Lip(s,L^q(S))$ is the Besov space $\mathrm{B}^s_{q,\infty}$~\cite[Ch.~2, Section~10]{devorelorentz93} which are very popular in approximation theory. In particular, $\Lip(s,L^{1/s}(S))$ contains the space $\BV(S)$ of functions of bounded variation on $S$; see~\cite[Ch.~2, Lemma~9.2]{devorelorentz93}. Thus Lipschitz spaces are rich enough to contain functions with both discontinuities and fractal structure.

We now state the following approximation error bounds whose proofs use standard arguments from approximation theory; see~\cite[Section~6.2.1]{Hafienne_FE_Nonlocal_plaplacian_evolu} for details.
\begin{lem} 
There exists a positive constant $C_s$, depending only on $s$, such that for all $F \in \Lip(s,L^q(S))$, $s \in ]0,1]$, $q \in [1,+\infty]$,
\begin{equation}
\anorm{F - \injn\projn F}_{L^q(S)} \le C_s \deltan^s \aabs{F}_{\Lip(s,L^q(S))}.
\label{eq:lipspaceapprox}
\end{equation}
\label{lem:spaceapprox}
\end{lem}

We denote by $\BV([0,T];L^q(S))$ the Banach space of functions $f: \O \times [0,T] \to \R$ such that
\[
\Var_{q}(f) \eqdef \sup_{0 \leq t_0 < t_1 < \cdots < t_N \leq T} \sum_{i=1}^N \norm{f(\cdot,t_i) - f(\cdot,t_{i-1})}_{L^q(S)} < +\infty, 
\]
endowed with the norm $\norm{f}_{\BV([0,T];L^q(S))} \eqdef \norm{f(0)}_{L^q(S)} + \Var_{q}(f)$.

\section{The sparse graph model}
\label{graphlimits}

One of our fundamental goals in this paper is to understand the behaviour of~\eqref{neumann} when discretized on a sequence of networks. Such networks can be of many types (biological, physical, social, data processing, etc.), whose details vary widely, but which bear similar structural phenomena. In turn, one may wonder whether discretization of~\eqref{neumann} is stable to the particular realization of the network and input data. It is then natural to consider a sequence of graphs with size tending to infinity and ask whether the solutions of discrete forms of~\eqref{neumann} on these graphs converge to any meaningful sort of limit. This will allow us in turn to establish a continuum limit of the solutions to these discrete evolution problems. For this, we need to be equipped with an appropriate theory of graph limits.

\subsection{$L^q$ graphons and graph limits}
\label{subsec:Lqgraphon}

In~\cite{Borgs18partII,Borgs19partI}, the authors laid the foundations of $L^q$ graphons, For $q \geq 1$, where the term $L^q$ graphon refers to a symmetric function $\K \in L^q([0,1]^2)$. They developed a theory of limits for sequences of sparse graphs based on such graphons, which generalizes both the existing theory of bounded graphons that are tailored to dense graph limits~\cite{LovaszBook}, and its extension in~\cite{Ballobas2008} to sparse graphs under a no dense spots assumptions. The latter graph model was studied in~\cite{HafieneRandom20} in the context of continuum limits of $p$-Laplacian evolution problems on graphs. Nevertheless, the boundedness assumption of the graphon underlying these graph models is still highly restrictive. In particular, it does not allow to handle singular graphons and corresponding network models which have statistics governed by power laws. The theory of unbounded by $L^q$ graphons allows to analyse graphs with power law degree distributions, hence providing a broadly applicable limit theory for sparse graphs with unbounded average degrees.

\subsection{Sparse $K$-random graph models}
\label{subsec:randomgraph}

We consider weighted graphs, which include as a special case simple unweighted graphs. Let $G = \pa{V(G),E(G)}$, be a weighted graphs with vertex set $V(G)$ and edge set $E(G) \subseteq V(G)^2$, respectively. In $G$, every edge $(i,j) \in E(G)$ (allowing loops with $i=j$) is given a weight $\beta_{ij} \in \R^+$\footnote{In \cite{Borgs19partI}, the weights are even allowed to be negative, but we will not consider this situation which is meaningless in our context.}. We set $\beta_{ij} = 0$ whenever $(i,j) \not\in E(G)$.

The theory of random graphs was founded in the 50's-60's by Erd\"os and R\'enyi~\cite{ErdosRenyi1960}, who started the systematic study of the space of graphs with $n$ labeled vertices and $M= M(n)$ edges, with all graphs equiprobable. The aim is to turn the set of all graphs with $n$ vertices into a probability space. Intuitively we should be able to generate a sequence of graphs $\seq{G_n}$ randomly as follows: for each edge $(i,j) \in [n]^2$, we decide by some random experiment whether or not $(i,j)$ shall be an edge of $G_n$, these experiments are performed independently.

The idea underlying the sparse $\K$-random graph model proposed by~\cite{Borgs19partI} is that each $L^q$ graphon $\K$ gives rise to a natural random graph model, which produces a sequence of sparse graphs converging to $\K$ in an appropriate metric. Inspired by their work, we propose the following construction. 

\begin{defi}
Fix $n \in \N^*$, let $\K$ be an $L^1$ graphon and $\rho_n > 0$. Take the equispaced partition of $[0,1]$ in intervals $]x_{i-1},x_i]$, $i \in [n]$, where $x_i = i/n$. Let $\bK \in \R_{+}^{n \times n}$ be a weight matrix such that:
\begin{enumerate}[label=\rm({\textbf{H${}_{w}$.\arabic*}})]
\item $\norm{\injn\bK - \K}_{L^1([0,1]^2)} \to 0$ as $n \to +\infty$. \label{assum:bKL1}
\item $\norm{\injn\bK(x,\cdot) - \K(x,\cdot)}_{L^1([0,1])} \to 0$ uniformly in $x \in [0,1]$. \label{assum:bKL1row}
\end{enumerate}

Generate the random graph 
\[
G_n = \pa{V(G_n),E(G_n)} \eqdef \bG(n,\K,\rho_n)
\] 
as follows: join each pair $(i,j) \in [n]^2$ of vertices independently, with probability 
\begin{equation}
\P\pa{(i,j) \in E(G_n) | \bX} = \rho_n \wedg_{ij} , \qwhereq \wedg_{ij} \eqdef \min \pa{\bK_{ij}, \rho_n^{-1}} .
\label{eq:sparsegraphmodelavg}
\end{equation}
\label{def:randomgraph}
\end{defi}

\begin{rem}
In the original sparse $\K$-random graph model defined in~\cite{Borgs19partI}, the $x_i$'s are random iid samples drawn from the uniform distribution on $[0,1]$. Moreover, $\bK_{ij} = \K(x_i,x_j)$. In this case, it follows from~\cite[Theorem~2.14(a)]{Borgs19partI} (which relies on~\cite[Theorem]{Hoeffding61}) that assumptions~\ref{assum:bKL1} holds with probability 1.

Another interesting case is where $\bK=\projn\K$. Thanks to Lemma~\ref{lem:projinj}, $\norm{\injn\projn\bK}_{L^1(\O^2)} \leq \norm{\K}_{L^1(\O^2)}$ with probability 1. Thus, the Lebesgue differentiation theorem and the dominated convergence theorem allow to assert that $\injn\projn\bK$ converges to $\K$ in $L^1(\O^2)$. In turn, assumption \ref{assum:bKL1} holds. 
\label{rem:graphmodel}
\end{rem}

For appropriate choices of $\rho_n$, the graph model constructed according to Definition~\ref{def:randomgraph} allows to sample both dense and sparse graphs from the graphon $\K$. In particular, the sparsity assumption $\rho_n \to 0$ reflects the fact that $\rho_n$ needs to be arbitrarily close to zero in order to see the unbounded/singular part of $\K$. The assumption that $n\rho_n \to +\infty$ means the average degree tends to infinity. To check this, the average number of edges in this graph model is
\begin{align*}
\EE\pa{E(\bG(n,\K,\rho_n))} 
&= \rho_n n^2 \pa{n^{-2} \sum_{(i,j) \in [n]^2} \wedg_{ij}} \\
&= \rho_n n^2 \pa{\norm{\injn\bK}_{L^1([0,1]^2)} - \norm{\pa{\injn\bK - \rho_n^{-1}}_+}_{L^1([0,1]^2)}}.
\end{align*}
By assumption~\ref{assum:bKL1}, we have $\norm{\injn\bK}_{L^1([0,1]^2)} = \norm{\K}_{L^1([0,1]^2)} + o(1)$. Moreover, since $\rho_n \to 0$, we have from~\eqref{eq:sumKrho} that $\norm{\pa{\injn\bK - \rho_n^{-1}}_+}_{L^1([0,1]^2)}=o(1)$. In turn,
\begin{align*}
\EE\pa{E(\bG(n,\K,\rho_n))} &= \rho_n n^2 \pa{\norm{\K}_{L^1([0,1]^2)} + o(1)} .
\end{align*}
As expected, this gives rise to a sparse graph whose edge density is $\rho_n \to 0$.
For the average degree of this graph model, arguing similarly to above, and using~\ref{assum:bKL1row}, we have
\begin{align*}
\EE\pa{\deg_{G_n}(i)} 
&= \rho_n n \pa{n^{-1} \sum_{j \in [n]} \wedg_{ij}} \\
&= \rho_n n \pa{\norm{\injn\bK(x_i,\cdot)}_{L^1([0,1])} - \norm{\pa{\injn\bK(x_i,\cdot) - \rho_n^{-1}}_+}_{L^1([0,1])}} \\
&= \rho_n n \pa{\int_{0}^1 \K(x_i,y) dy + o(1)} .
\end{align*}
As anticipated, the average degree is indeed unbounded since $\rho_n n \to +\infty$ .\\


The above sequence of graphs generated also enjoys the following convergence result.
\begin{prop}\label{prop:sparsegraphconv}
Let $\K$ be an $L^1$ graphon and $\bK$ be a weight matrix such that~\ref{assum:bKL1} holds. If $\rho_n > 0$ with $\rho_n \to 0$ and $n\rho_n \to +\infty$ as $n \to +\infty$, then $\rho_n^{-1}\bG(n,\K,\rho_n)$ converges almost surely to $\K$ in the cut distance metric (see~\cite{Borgs19partI,Borgs18partII} for details about this metric). 
\end{prop}
\bpf{}
We essentially adapt the arguments of in the proof of~\cite[Theorem~2.14(b)]{Borgs19partI}. More precisely, since~\ref{assum:bKL1} holds, one has to show~\cite[(7.1)]{Borgs19partI}. For this, we invoke~\cite[Lemma~7.3]{Borgs19partI} by checking the condition (7.3) therein. We have by sublinearity of $(\cdot)_+$ that
\begin{equation}\label{eq:sumKrho}
\begin{aligned}
\frac{1}{n^2} \sum_{(i,j) \in [n]^2} \pa{\bK_{ij} - \rho_n^{-1}}_+ 
&=\int_{[0,1]^2} \pa{\injn\bK(x,y) - \rho_n^{-1}}_+ dxdy \\
&\leq \int_{[0,1]^2} \pa{\injn\bK(x,y) - \K(x,y)}_+ dxdy + \int_{[0,1]^2} \pa{\K(x,y) - \rho_n^{-1}}_+ dxdy \\
&\leq \norm{\injn\bK - \K}_{L^1([0,1]^2)} + \int_{[0,1]^2} \pa{\K(x,y) - \rho_n^{-1}}_+ dxdy .
\end{aligned}
\end{equation}
The right-hand side in the above display goes to $0$ as $n \to +\infty$ by \ref{assum:bKL1} and since $\rho_n \to 0$. Indeed, for every $L > 0$, the limit superior of the last term is bounded by $\norm{(\K - L)_+}_{L^1([0,1]^2})$, and this can be made arbitrarily small by choosing $L$ large.
\epf{}

\begin{exe}
For an example that cannot be handled using $L^\infty$ graphons, and thus does not enter in the framework of \cite{Hafienne_FE_Nonlocal_plaplacian_evolu,HafieneRandom20}, consider a $\K$-random graph model $\bG(n,\K,\rho_n)$ constructed according to Definition~\ref{def:randomgraph} with $\bK=\projn\K$, where $\K(x,y) = \J(x-y)$, $\J: z \in [-1,1] \mapsto 2^{-1}(1-\beta)(2-\beta) |z|^{-\beta}$, $\beta \in ]0,1[$. First, observe that the radially symmetric kernel $\J$ is singular but fulfills all assumptions \ref{assum:Jpos}, \ref{assum:Jsym} and \ref{assum:JL1}. In addition, by virtue of Remark~\ref{rem:graphmodel}, \ref{assum:bKL1}-\ref{assum:bKL1row} also hold with
\[
\norm{\K}_{L^1([0,1]^2)} = 1 \qandq \int_{0}^1 \K(x,y) dy = 2^{-1}(2-\beta)\pa{x^{1-\beta} + (1-x)^{1-\beta}} \in 2^{-1}(2-\beta)[1,2^\beta].
\]
\end{exe}

We also have the following convergence result in the $L^{\infty,1}$ norm that will be instrumental in Section~\ref{sec:applicationgraphs}. According to the construction in Definition~\ref{def:randomgraph}, we let $\bLam_{ij}$, $(i,j) \in [n]^2, i \neq j$, be random variables such that $\rho_n \bLam_{ij}$ follows a Bernoulli distribution with parameter $\rho_n \wedg_{ij}$. For each row $i \in [n]$, $\pa{\bLam_{ij}}_{j \in [n]}$ are independent.
\begin{lem}
Let $\K$ be an $L^{\infty,1}$ graphon, i.e. it satisfies~\ref{assum:Kpos},~\ref{assum:Ksym} and~\ref{assum:KL1}. Take the weight matrix $\bK=\projn\K$. Assume that $\rho_n \to 0$ and $n\rho_n = \omega\pa{\pa{\log n}^\gamma}$ for some $\gamma > 1$. Then with probability 1,
\[
\norm{\injn\bLam}_{L^{\infty,1}([0,1]^2)} - \norm{\injn\wedg}_{L^{\infty,1}([0,1]^2)} \to 0 .
\]
If, moreover, \ref{assum:bKL1row} holds, then
\[
\norm{\injn\bLam}_{L^{\infty,1}([0,1]^2)} \to \norm{\K}_{L^{\infty,1}([0,1]^2)} .
\]
with probability 1.
\label{lem:Linf1norm}
\end{lem}

\bpf{}
For any $\varepsilon > 0$, we have by the union bound
\begin{align*}
&\P\pa{\abs{\norm{\injn\bLam}_{L^{\infty,1}([0,1]^2)} - \norm{\injn\wedg}_{L^{\infty,1}([0,1]^2)}} \geq \varepsilon} \\
&=\P\pa{\abs{\max_{i}\sum_{j}\bLam_{ij} - \max_{i}\sum_{j}\wedg_{ij}} \geq \varepsilon n} \\
&=\P\pa{\abs{\max_{i}\sum_{j}\rho_n\bLam_{ij} - \max_{i}\sum_{j}\rho_n\wedg_{ij}} \geq \varepsilon\rho_n n} \\
&\leq \P\pa{\max_{i}\abs{\sum_{j}\rho_n(\bLam_{ij} - \wedg_{ij})} \geq \varepsilon\rho_n n} \\
&\leq \sum_i\P\pa{\abs{\sum_{j}\rho_n(\bLam_{ij} - \wedg_{ij})} \geq \varepsilon\rho_n n} .
\end{align*}
Since $\pa{\rho_n\bLam_{ij}}_{j}$ are independent Bernoulli variables with means $\pa{\rho_n \wedg_{ij}}_j$, it follows from the variant of the Chernoff bound in~\cite[Lemma~7.1]{Borgs19partI}, that for every $\varepsilon > 0$,
\begin{align*}
&\P\pa{\abs{\norm{\injn\bLam}_{L^{\infty,1}([0,1]^2)} - \norm{\injn\wedg}_{L^{\infty,1}([0,1]^2)}} \geq \varepsilon} \\
&\leq 2\sum_i\exp\pa{-\frac{1}{3}\min\pa{\frac{\varepsilon\rho_nn}{\rho_n\sum_{j}\wedg_{ij}},1}\varepsilon\rho_n n} \\
&\leq 2n\exp\pa{-\frac{1}{3}\min\pa{\frac{\varepsilon}{\norm{\injn\wedg}_{L^{\infty,1}([0,1]^2)}},1}\varepsilon\rho_n n} \\
&\leq 2n\exp\pa{-\frac{1}{3}\min\pa{\frac{\varepsilon}{\norm{\K}_{L^{\infty,1}([0,1]^2)}},1}\varepsilon \omega\pa{(\log n)^\gamma}} \\
&\leq 2n^{-\omega\pa{(\log n)^{\gamma-1}}} ,
\end{align*}
since $\gamma > 1$, and where we used~\eqref{eq:sparsegraphmodelavg} and Lemma~\ref{lem:projinj} to show that
\[
\norm{\injn\wedg}_{L^{\infty,1}([0,1]^2)} \leq \norm{\injn\bK}_{L^{\infty,1}([0,1]^2)} = \norm{\injn\projn\K}_{L^{\infty,1}([0,1]^2)} \leq \norm{\K}_{L^{\infty,1}([0,1]^2)} .
\]
Invoking the (first) Borel-Cantelli lemma, we have the first claim. On the other hand,
\begin{align*}
&\abs{\norm{\injn\wedg}_{L^{\infty,1}([0,1]^2)} - \norm{\K}_{L^{\infty,1}([0,1]^2)}} 
\leq \norm{\injn\wedg - \K}_{L^{\infty,1}([0,1]^2)} \\
&\leq \norm{\injn\wedg - \injn\projn\K}_{L^{\infty,1}([0,1]^2)} + \norm{\injn\projn\K - \K}_{L^{\infty,1}([0,1]^2)} \\
&= \norm{(\injn\projn\K - \rho_n^{-1})_+}_{L^{\infty,1}([0,1]^2)} + \norm{\injn\projn\K - \K}_{L^{\infty,1}([0,1]^2)} \\
&\leq \norm{(\K - \rho_n^{-1})_+}_{L^{\infty,1}([0,1]^2)} + \norm{(\injn\projn\K - \K)_+}_{L^{\infty,1}([0,1]^2)} + \norm{\injn\projn\K - \K}_{L^{\infty,1}([0,1]^2)} \\
&\leq \norm{(\K - \rho_n^{-1})_+}_{L^{\infty,1}([0,1]^2)} + 2\norm{\injn\projn\K - \K}_{L^{\infty,1}([0,1]^2)} .
\end{align*}
Since $\rho_n \to 0$ and in view of \ref{assum:bKL1row}, the right-hand side in the above display goes to $0$ as $n \to +\infty$. Combined with the first claim we obtain the desired conclusion.
\epf{}

\section{Well-posedness}
\label{existenceuniqueness}

\subsection{The case $p \in ]1,+\infty[$}
To lighten notation, for $ 1< p < +\infty$, we define the function 
\begin{equation*}
\begin{split}
\Psi : x \in \mathbb{R} \mapsto \abs{x}^{p-2} x = \sign(x) {\abs{x}}^{p-1} ,
\end{split}
\end{equation*}
where we take $\sign(0)=0$.

The next lemma summarizes key monotonicity and continuity properties of $\Psi$ which will be instrumental to us.
\begin{lem}\label{lem:psiprop}
\begin{enumerate}[label={\rm (\roman{*})}]
\item Monotonicity: assume that the constant $\beta$ satisfies $\beta\in [\max(p,2),+\infty[$. Then for all $x,y \in \R$, \label{item:psimonotone}
\begin{equation}
\pa{\Psi(y)-\Psi(x)}(y-x) \geq C_1
\abs{y-x}^\beta\pa{|y|+|x|}^{p-\beta},
\label{eq:psigenmonotone}
\end{equation}
where the constant $C_1$ is sharp and given by
\begin{equation}
C_1= 2^{2-p}\min (1,p-1).
\label{eq:psigenmonotonecst}
\end{equation}
In particular,
\begin{equation}\label{eq:psimonotone}
\pa{\Psi(y)-\Psi(x)}(y-x) \geq C_1
\begin{cases}
\abs{y-x}^p 			& p \in [2,+\infty[ , \\
\abs{y-x}^2\pa{|y|+|x|}^{p-2} 	& p \in ]1,2].
\end{cases}
\end{equation}

%

\item Continuity: assume that the constant $\alpha$ satisfies $\alpha\in [0,\min(1,p-1)]$. Then for all $x,y \in \R$,  \label{item:psicont}
\begin{equation}
\abs{\Psi(y)-\Psi(x)} \leq C_2
\abs{y-x}^\alpha \pa{|y|+|x|}^{p-1-\alpha},
\label{eq:psigencont}
\end{equation}
where the constant $C_2$ is sharp and given by
\begin{equation}
C_2= \max (2^{2-p},(p-1) 2^{2-p},1).
\label{eq:psigencontcst}
\end{equation}
In particular,
\begin{equation}\label{eq:psicont}
\abs{\Psi(y)-\Psi(x)} \leq C_2
\begin{cases}
\abs{y-x} \pa{|y|+|x|}^{p-2} 	& p \in [2,+\infty[ , \\
\abs{y-x}^{p-1} 		& p \in ]1,2] ,
\end{cases}
\end{equation}
\end{enumerate}
\end{lem}

\bpf{}
\begin{enumerate}[label={\rm (\roman{*})}]
\item For \eqref{eq:psigenmonotone}, see \cite[Theorem~2.2]{bystrom2005sharp}. For~\eqref{eq:psimonotone}, set $\beta=p$ for $p \geq 2$ and $\beta=2$ otherwise in \eqref{eq:psigenmonotone}; see also the seminal results of~\cite[Lemma~5.1 and Lemma~5.2]{Glowinski75}.
\item For \eqref{eq:psigencont}, see \cite[Theorem~2.1]{bystrom2005sharp}. For~\eqref{eq:psicont}, set $\alpha=1$ for $p \geq 2$ and $\alpha=p-1$ otherwise in \eqref{eq:psigencont}; see also the seminal results of~\cite[Lemma~5.3 and Lemma~5.4]{Glowinski75}.
\end{enumerate}
\epf{}


We now collect some preliminary properties of the nonlocal $p$-Laplacian, an operator on $L^1(\O)$ that we denote for short as
\[
\plap{\K}{p} u(\bx,t) = - \int_{\O} \K(\bx,\by) \abs{u(\by,t) - u(\bx,t) }^{p-2} (u(\by,t) - u(\bx,t) ) d\by .
\]

\begin{prop}
Assume that $\K$ satisfies~\ref{assum:Kpos},~\ref{assum:Ksym} and~\ref{assum:KL1}.
\begin{enumerate}[label={\rm (\roman{*})}]
\item $\plap{\K}{p}$ is positively homogeneous of degree $p-1$. \label{item1homo}
\item If $p > 2$, $L^{p-1} (\O) \subset \dom(\plap{\K}{p})$ . \label{item2dom} 
\item If $1< p \leq 2 $, $\dom(\plap{\K}{p}) = L^1(\O)$ and $\plap{\K}{p}$ is closed in $L^1(\O) \times L^1(\O)$. \label{item3dom}
\item Let $h: \R \to \R$. Then for every $u,v \in L^p(\O)$, \label{itemmonotonep}
\begin{equation}
\hspace*{-1cm}
\begin{split}
0 &\leq \int_{\O} \pa{\plap{\K}{p} u(\bx) - \plap{\K}{p} v(\bx)}h(u(\bx) - v(\bx)) d\bx \\
 &= \frac{1}{2} \int_{\O^2} \K(\bx,\by) \pa{\F(u(\by) - u(\bx)) - \F(v(\by) - v(\bx))} \pa{h(u(\by) - v(\by)) - h(u(\bx) - v(\bx))} d\by d\bx .
\end{split}
\label{eq:monotoneplap}
\end{equation}
If $h$ is bounded, then this holds for any $u,v \in \dom(\plap{\K}{p})$.

\item For every $u,v \in L^p(\O)$, \label{itemstrongmonotonep}
\begin{multline*}
\int_{\O} \pa{\plap{\K}{p} u(\bx) - \plap{\K}{p} v(\bx))(u(\bx) - v(\bx))} d\bx \geq \\
\frac{C}{2} \pa{\int_{\O^2} \K(\bx,\by) \abs{(u(\by) - u(\bx)) - (v(\by) - v(\bx))}^p d\by d\bx}^{\max(1,2/p)}
\end{multline*}
where
\[
C = 
\begin{cases}
C_1 & p \in [2,+\infty[ , \\
2^{2p-5}C_1 \norm{\K}_{L^{\infty,1}(\O^2)}^{1-2/p}\pa{\norm{u}_{L^p(\O)}+\norm{v}_{L^p(\O)}}^{p-2} & p \in ]1,2[ . 
\end{cases}
\]
and $C_1$ is the constant in \eqref{eq:psigenmonotonecst}. If $u,v \in L^\infty(\O)$, then
\begin{multline*}
\int_{\O} \pa{\plap{\K}{p} u(\bx) - \plap{\K}{p} v(\bx))(u(\bx) - v(\bx))} d\bx \geq \\
\frac{C}{2} \pa{\int_{\O^2} \K(\bx,\by) \abs{(u(\by) - u(\bx)) - (v(\by) - v(\bx))}^2 d\by d\bx}^{\max(1,p/2)},
\end{multline*}
where
\[
C = 
\begin{cases}
C_1 & p \in [2,+\infty[ , \\
2^{p-2}C_1 \pa{\norm{u}_{L^\infty(\O)}+\norm{v}_{L^\infty(\O)}}^{p-2} & p \in ]1,2[ . 
\end{cases}
\]

\item For $p \in ]1,2]$ and every $u,v \in L^2(\O)$, \label{itemholderp}
\begin{align*}
\int_{\O} \pa{\plap{\K}{p} u(\bx) - \plap{\K}{p} v(\bx))(u(\bx) - v(\bx))} d\bx \geq
C \norm{\plap{\K}{p} u - \plap{\K}{p} v}_{L^2(\O)}^{p/(p-1)}
\end{align*}
where
\[
C = 2^{\frac{p-2}{2(p-1)}} \pa{C_2^{1/2}\norm{\K}_{L^{\infty,1}(\O^2)}}^{\frac{1}{1-p}}(1-1/p), \qandq \text{ $C_2$ is the constant in \eqref{eq:psigencontcst}} .
\]

\item For $p \in ]1, +\infty[$, $\plap{\K}{p}$ is completely accretive and satisfies the range condition \label{item4compaccretive}
\begin{equation}
L^p(\O) \subset \range(\Id + \plap{\K}{p}).
\label{eq:rangecondition}
\end{equation}
Consequently, the resolvent $J_{\lambda\plap{\K}{p}} \eqdef \pa{\Id+\lambda\plap{\K}{p}}^{-1}$, $\lambda > 0$, is single-valued on $L^p(\O)$ and nonexpansive in $L^q (\O)$ for all $q \in [1,+\infty]$. 
\end{enumerate}
\label{prop:plap}
\end{prop}

\bpf{}
\ref{item1homo},~\ref{item2dom} and~\ref{item3dom} follow from~\cite[Remark~2.2]{mazonjmaa} which still holds for our larger class of kernels $\K$.

For~\ref{itemmonotonep}, see~\cite[Lemma~A.2]{Hafienne_FE_Nonlocal_plaplacian_evolu}. Monotonicity is immediate since $h$ is non-decreasing.

The proof of~\ref{item4compaccretive} is the same as that of~\cite[Theorem~2.4]{mazonjmaa}, where we invoke the monotonicity claim~\ref{item1homo}. 

We now show \ref{itemstrongmonotonep}\footnote{This can be seen as a nonlocal analogue of~\cite[Proposition~5.1 and Proposition~5.2]{Glowinski75}.}. The case $p \in [2,+\infty[$ is immediate by inserting Lemma~\ref{lem:psiprop}\ref{item:psimonotone} into \eqref{eq:monotoneplap} with $h(x)=x$. For $p \in ]1,2]$, to lighten notation, denote the nonlocal gradient $\gradNL u(\bx,\by)=u(\by)-u(\bx)$. We then have by Lemma~\ref{lem:psiprop}\ref{item:psimonotone} that
\begin{multline}\label{eq:gradmonotonep}
C_1|\gradNL (u-v)(\bx,\by)|^2 \leq \\
\pa{\Psi(\gradNL u(\bx,\by))-\Psi(\gradNL v(\bx,\by))}(\gradNL u(\bx,\by)-\gradNL v(\bx,\by)) \pa{|\gradNL u(\bx,\by)|+|\gradNL v(\bx,\by)|}^{2-p}.
\end{multline}
Taking the power $p/2$, multiplying by $\K$ and integrating, we get
\begin{multline*}
C_1^{p/2} \int_{\O^2} \K(\bx,\by)|\gradNL (u-v)(\bx,\by)|^p d\bx d\by
\leq \\
\int_{\O^2} \pa{\K(\bx,\by)\pa{\Psi(\gradNL u(\bx,\by))-\Psi(\gradNL v(\bx,\by))}(\gradNL u(\bx,\by)-\gradNL v(\bx,\by))}^{p/2} \\
\pa{\K(\bx,\by)^{1/p}(|\gradNL u(\bx,\by)|+|\gradNL v(\bx,\by)|)}^{(2-p)p/2} d\bx d\by .
\end{multline*}
It is easily seen that
\begin{align*}
\pa{\K \cdot \pa{\Psi(\gradNL u)-\Psi(\gradNL v)}(\gradNL u-\gradNL v)}^{p/2} &\in L^{2/p}(\O^2) \\
\pa{\K^{1/p} \cdot (|\gradNL u|+|\gradNL v|)}^{(2-p)p/2} &\in L^{2/(2-p)}(\O^2) .
\end{align*}
It then follows from H\"older inequality and~\eqref{eq:monotoneplap} that
\begin{multline*}
C_1^{p/2} \int_{\O^2} \K(\bx,\by)|\gradNL (u-v)(\bx,\by)|^p d\bx d\by
\leq \\
2\pa{\int_{\O}\pa{\plap{\K}{p} u(\bx) - \plap{\K}{p} v(\bx)}(u(\bx) - v(\bx)) d\bx}^{p/2} \cdot \\
\pa{\int_{\O^2}\K(\bx,\by)(|\gradNL u(\bx,\by)|+|\gradNL v(\bx,\by)|)^p d\bx d\by}^{(2-p)/2} .
\end{multline*}
We have by Jensen's inequality
\begin{align*}
&\int_{\O^2}\K(\bx,\by)(|\gradNL u(\bx,\by)|+|\gradNL v(\bx,\by)|)^p d\bx d\by \\
&\leq 4^{p-1}\int_{\O^2}\K(\bx,\by)(|u(\bx)|^p+|u(\by)|^p+|v(\bx)|^p+|v(\by)|^p) d\bx d\by \\
&\leq 2^{2p-1}\norm{\K}_{L^{\infty,1}(\O^2)}\pa{\norm{u}_{L^p(\O)}^p+\norm{v}_{L^p(\O)}^p} ,
\end{align*}
whence we obtain
\begin{multline*}
C_1 \pa{\int_{\O^2} \K(\bx,\by)|\gradNL (u-v)(\bx,\by)|^p d\bx d\by}^{2/p}
\leq \\
2^{5-2p}\pa{\int_{\O}\pa{\plap{\K}{p} u(\bx) - \plap{\K}{p} v(\bx)}(u(\bx) - v(\bx)) d\bx} \\
\norm{\K}_{L^{\infty,1}(\O^2)}^{(2-p)/p}\pa{\norm{u}_{L^p(\O)}+\norm{v}_{L^p(\O)}}^{2-p} .
\end{multline*}
Rearranging proves the bound. For $u,v \in L^\infty(\O)$ and $p \in [2,+\infty]$ we use that $L^p(\O) \subset L^2(\O)$. For $p \in ]1,2]$, we embark from~\eqref{eq:gradmonotonep} and use that for all $(\bx,\by) \in \O^2$,
\begin{align*}
|\gradNL u(\bx,\by)|+|\gradNL v(\bx,\by)| \leq 2\pa{\norm{u}_{L^\infty(\O)}+\norm{v}_{L^\infty(\O)}} .
\end{align*}
Multiplying~\eqref{eq:gradmonotonep} by $\K$, integrating and using~\eqref{eq:monotoneplap}, we conclude.

To prove~\ref{itemholderp}, we start by showing that $\plap{\K}{p}$ is H\"older continuous with exponent $p-1$ on $L^2(\O)$. We have by Jensen inequality (twice) and \eqref{eq:psicont},
\begin{align}
&\norm{\plap{\K}{p} u - \plap{\K}{p} v}_{L^2(\O)}^2 
= \int_{\O} \aabs{\int_{\O} \K(\bx,\by)\pa{\Psi(\gradNL u(\bx,\by)) - \Psi(\gradNL v(\bx,\by))} d\by}^2 d\bx \nonumber\\
&\leq \norm{\K}_{L^{\infty,1}(\O^2)}\int_{\O^2}\K(\bx,\by)\pa{\Psi(\gradNL u(\bx,\by)) - \Psi(\gradNL v(\bx,\by))}^2 d\bx d\by \nonumber\\
&\leq C_2\norm{\K}_{L^{\infty,1}(\O^2)}\int_{\O^2}\K(\bx,\by)\pa{\gradNL (u-v)(\bx,\by))}^{2(p-1)} d\bx d\by \nonumber\\
&\leq 2^p C_2\norm{\K}_{L^{\infty,1}(\O^2)}\int_{\O^2}\K(\bx,\by)\pa{u(\bx)-v(\bx)}^{2(p-1)} d\bx d\by \nonumber\\
&\leq 2^p C_2\norm{\K}_{L^{\infty,1}(\O^2)}^2\pa{\int_{\O}\pa{u(\bx)-v(\bx)}^{2} d\bx}^{p-1} \nonumber\\
&= 2^p C_2\norm{\K}_{L^{\infty,1}(\O^2)}^2\norm{u-v}_{L^2(\O)}^{2(p-1)} . \label{eq:plapholder}
\end{align}
We are now in position to invoke~\cite[Corollary~18.14(i)$\Rightarrow$(v)]{bauschke2011convex} to show that the claimed inequality holds.
\epf{}

Solutions of~\eqref{neumann} will be understood in the following sense: 
\begin{defi}\label{def:sol}
Let $p \in ]1,+\infty[$. A solution of~\eqref{neumann} in $[0,T]$ is a function 
\[
u \in C([0,T];L^1(\O)) \cap W^{1,1}(]0,T[;L^1(\O)),
\]
that satisfies $u(\bx,0) = g(\bx)$ a.e. $\bx \in \O$ and 
\[
\frac{\partial }{\partial t} u(\bx,t) = -\plap{\K}{p} u(\bx,t) + f(\bx,t) \quad \text{ a.e. in } \O \times ]0,T[ . 
\]
Such a solution is also a strong solution (see~\cite[Definition~A.3]{vaillo}).
\end{defi}

The main result of existence and uniqueness of a global solution, that is, a solution on $[0,T]$ for $T>0$ is stated in the following theorem.

\begin{theo}
\label{theo:exist}
Suppose that $p \in ]1,+\infty[$ and assumptions~\ref{assum:Kpos},~\ref{assum:Ksym} and~\ref{assum:KL1} hold. Let $g \in L^p(\O)$ and $f \in L^1([0,T];L^p(\O))$.
\begin{enumerate}[label=(\roman*)]
\item For any $T > 0$, there exists a unique strong solution in $[0,T]$ of~\eqref{neumann}.
\item Moreover, for $q \in [1,+\infty]$, if $g_i \in L^q(\O)$ and $f_i \in L^1([0,T];L^q(\O))$, $i=1,2$, and $u_i$ is the solution of \eqref{neumann} with data $(f_i,g_i)$, then
\begin{equation}
\norm{u_1(\cdot,t) - u_2(\cdot,t)}_{L^q(\O)} \le \norm{g_1 - g_2}_{L^q(\O)} + \norm{f_1 - f_2}_{L^1([0,T];L^q(\O))}, \quad \forall t \in [0,T] .
\label{contraction}
\end{equation}
\end{enumerate}
\label{existence}
\end{theo}

\bpf{}
The proof follows the same lines as that of~\cite[Theorem~1.2]{mazonjmaa} extended to the case where $f \not\equiv 0$ thanks to the results of~\cite{Benilan72}, where we invoke Proposition~\ref{prop:plap}\ref{item2dom},~\ref{item3dom} and~\ref{item4compaccretive}.
\epf{}

\begin{rem}
In~\cite{mazonjmaa} (see also~\cite[Chapte~6]{vaillo}), the authors impose the following stringent assumptions: $\K(\bx,\by)=\J(\bx-\by)$, where $\J$ is nonnegative, continuous, radially symmetric, compactly supported, $\J(0) > 0$ and $\int_{\R^d} \J(\bx) d\bx < +\infty$. Actually, these assumptions are not needed for existence and uniqueness. The particular form $\J(\bx-\by)$ of the kernel is not needed. Continuity with radial symmetry and support compactness play a pivotal role to study convergence to the local $p$-Laplacian problem in~\cite[Theorem~1.5]{mazonjmaa}. In addition, $\J(0) > 0$ was mandatory to prove a Poincar\'e-type inequality in~\cite[Proposition~4.1]{mazonjmaa}. Even for the form $\J(\bx-\by)$, our assumptions~\ref{assum:Jpos},~\ref{assum:Jsym} and~\ref{assum:JL1} are weaker than those of~\cite{mazonjmaa}. This discussion remains true also for the case $p=1$.
\end{rem}

\subsection{The case $p=1$}
We will need to define subdifferential of the absolute value function on $\R$, which is the well-known set-valued mapping $\partial |\cdot|: \R \to 2^\R$,
\[
\partial |\cdot|(x) = 
\begin{cases}
1 	& x > 0 \\
[-1,1] 	& x = 0 \\
-1	& x < 0 .
\end{cases}
\]
It will be convenient to denote the $1$-Laplacian $\plap{\K}{1}$. This is a set-valued operator in $L^1(\O) \times L^1(\O)$ such that $\eta \in \plap{\K}{1} u$ if and only if
\[
\eta(\bx) = -\int_{\O} \K(\bx,\by) w(\bx,\by) d\by \quad \text{ a.e. in } \O, 
\]
for a subgradient function $w$ satisfying $\norm{w}_{L^\infty(\O^2)} \leq 1$, $w(\bx,\by)=-w(\by,\bx)$, and 
\[
w(\bx,\by) \in \partial |\cdot| (u(\by)-u(\bx)) .
\]

Solutions of~\eqref{neumann} will be understood in the following sense. 
\begin{defi}\label{def:solp1}
A solution of~\eqref{neumann} for $p=1$ in $[0,T]$ is a function 
\[
u \in C([0,T];L^1(\O)) \cap W^{1,1}(]0,T[;L^1(\O)),
\]
that satisfies $u(\bx,0) = g(\bx)$ for a.e. $\bx \in \O$ and 
\[
\frac{\partial }{\partial t} u(\bx,t) = -\eta(\bx,t) + f(\bx,t) \quad \text{ a.e. in } \O \times ]0,T[ , 
\]
where $\eta(\cdot,t) \in \plap{\K}{1} u(\cdot,t)$.
\end{defi}

Observe that for $p=1$, the evolution problem~\eqref{neumann} reads
\begin{equation*}
\begin{cases}
\frac{\partial }{\partial t} u(\bx,t) = \int_{\O} \K(\bx,\by) \sign(u(\by,t) - u(\bx,t)) d\by + f(\bx,t) , & \bx \in \O, t > 0,
\\
u(\bx,0)= g(\bx), \quad \bx \in \O,
\end{cases}
\end{equation*}
where 
\[
\sign(x) = 
\begin{cases}
\frac{x}{|x|} 	& x \neq 0 \\
0	& x = 0 .
\end{cases}
\]
Thus, it satisfies 
\[
\frac{\partial }{\partial t} u(\cdot,t) \in -\plap{\K}{1} u(\cdot,t) .
\]

In the same vein as Proposition~\ref{prop:plap}, the $1$-Laplacian enjoys the following properties.

\begin{prop}
Assume that $\K$ satisfies \ref{assum:Kpos}, \ref{assum:Ksym} and \ref{assum:KL1}.
\begin{enumerate}[label={\rm (\roman{*})}]
\item $\dom(\plap{\K}{1}) = L^1(\O)$ and (the graph of) $\plap{\K}{1}$ is closed in $L^1(\O) \times L^1(\O)$. \label{itemdomp1}
\item Let $h \in C^1(\R)$ be a nondecreasing function. Then for every $u_i \in L^1(\O)$ and any $\eta_i \in \plap{\K}{1}u_i$, $i=1,2$, \label{itemmonotonep1}
\begin{equation}
\begin{split}
0 &\leq \int_{\O} \pa{\eta_1(\bx) - \eta_2(\bx)}\pa{h(u_1(\bx) - u_2(\bx))} d\bx \\
 &= \frac{1}{2} \int_{\O^2}\K(\bx,\by) \pa{w_1(\bx,\by) - w_2(\bx,\by)}\pa{h(u_1(\by) - u_2(\by)) - h(u_1(\bx) - u_2(\bx))} d\bx d\by .
\end{split}
\label{eq:monotonep1}
\end{equation}
where $w_i$ are the corresponding subgradient functions defined above. 
In particular,
\[
\int_{\O^2} \K(\bx,\by) w_i(\bx,\by) u_i(\bx) d\bx d\by = -\frac{1}{2} \int_{\O^2} \K(\bx,\by) \abs{u_i(\by)-u_i(\bx)} d\bx d\by .
\]
\item $\plap{\K}{1}$ is completely accretive and satisfies the range condition \label{itemcompaccretivep1}
\begin{equation}
L^\infty(\O) \subset \range(\Id + \plap{\K}{1}).
\label{eq:rangeconditionp1}
\end{equation}
\end{enumerate}
\label{prop:1lap}
\end{prop}

\bpf{}
For~\ref{itemdomp1}, see~\cite[Remark~2.8]{mazonjmaa} which still holds for our class of kernels $\K$.

The proof of~\ref{itemcompaccretivep1} is again the same as that of~\cite[Theorem~2.9]{mazonjmaa}, where we invoke the monotonicity claim~\ref{itemmonotonep1} to which we turn now.

For any $v \in L^1(\O)$, we have the integration by parts formula
\begin{align}
&\int_{\O^2} \K(\bx,\by) w_i(\bx,\by) (v(\by)-v(\bx)) d\bx d\by \\
&=-\int_{\O^2} \K(\by,\bx) w_i(\by,\bx) v(\by)d\by d\bx - \int_{\O^2} \K(\bx,\by) w_i(\bx,\by)v(\bx)) d\bx d\by \\
&=-2\int_{\O^2} \K(\bx,\by) w_i(\bx,\by)v(\bx) d\bx d\by .
\label{eq:integpartsp1}
\end{align}
Taking $v(\bx)=h(u_1(\bx) - u_2(\bx))$ in~\eqref{eq:integpartsp1} with $w_1$ and $w_2$, and then taking the difference, we arrive at
\begin{align*}
&-2\int_{\O} \pa{\int_{\O}\K(\bx,\by) (w_1(\bx,\by)-w_2(\bx,\by))d\by} h(u_1(\bx) - u_2(\bx)) d\bx \\
&=2 \int_{\O} \pa{\eta_1(\bx) - \eta_2(\bx)}\pa{h(u_1(\bx) - u_2(\bx))} d\bx \\
&= \int_{\O^2} \K(\bx,\by) \pa{w_1(\bx,\by)-w_2(\bx,\by)}\pa{h(u_1(\by) - u_2(\by))-h(u_1(\bx) - u_2(\bx))} d\bx d\by .
\end{align*}
By the mean-value theorem applied to $h$, we get
\begin{align*}
&=2 \int_{\O} \pa{\eta_1(\bx) - \eta_2(\bx)}\pa{h(u_1(\bx) - u_2(\bx))} d\bx \\
&= \int_{\O^2} \K(\bx,\by) \pa{w_1(\bx,\by)-w_2(\bx,\by)}h'(\zeta(\bx,\by))\pa{(u_1(\by) - u_2(\by))-(u_1(\bx) - u_2(\bx))} d\bx d\by \\
&= \int_{\O^2} \K(\bx,\by) h'(\zeta(\bx,\by)) \pa{w_1(\bx,\by)-w_2(\bx,\by)}\pa{(u_1(\by) - u_1(\bx)) - (u_2(\by)-u_2(\bx))} d\bx d\by ,
\end{align*}
where $\zeta(\bx,\by)$ is an intermediate value between $u_1(\by) - u_2(\by)$ and $u_1(\bx) - u_2(\bx)$. Since $h$ is increasing, that $w_i(\bx,\by) \in \partial |\cdot| (u_i(\by)-u_i(\bx))$, and $\partial |\cdot|$ is a monotone operator, we get the claimed monotonicity.

To get the particular identity, we specialize~\eqref{eq:integpartsp1} by taking $v=u_i$, which entails
\begin{align*}
-\int_{\O^2} \K(\bx,\by) w_i(\bx,\by) (u_i(\by)-u_i(\bx)) d\bx d\by 
=2\int_{\O^2} \K(\bx,\by) w_i(\bx,\by)u_i(\bx) d\bx d\by .
\end{align*}
We finally use the equivalent characterization of $\partial \abs{\cdot}$, which originates from the Fenchel's identity since $\abs{\cdot}$ is positively homogeneous,
\[
\partial \abs{\cdot}(x) = \ens{\xi \in \R}{\abs{\xi} \leq 1 \qandq \xi x = \abs{x}} .
\]
Applying this identity with $x=u_i(\by)-u_i(\bx)$ and $\xi=w_i(\bx,\by)$ gives the claim.
\epf

\begin{theo}
\label{theo:existp1}
Suppose that $p=1$, and assumptions~\ref{assum:Kpos},~\ref{assum:Ksym} and~\ref{assum:KL1} hold. Let $g \in L^1(\O)$ and $f \in L^1([0,T];L^1(\O))$. For any $T > 0$, there exists a unique solution in $[0,T]$ of~\eqref{neumann} in the sense of Definition~\ref{def:solp1}.
\label{existencep1}
\end{theo}

\bpf{}
The proof is an adaptation of~\cite[Theorem~1.4]{mazonjmaa} to the case where $f \not\equiv 0$ thanks to the results of~\cite{Benilan72}, where we invoke Proposition~\ref{prop:1lap}\ref{itemdomp1} and~\ref{itemcompaccretivep1}.
\epf{}

\section{Continuous-continuous estimates}
\label{sec:consistencecontinuous}

In this section, we provide an estimate that compares solutions of two $p$-Laplacian evolution problems of the form~\eqref{neumann} with two different kernels and initial data. This estimate will be instrumental to derive error bounds in the totally discrete case.

\subsection{The case $p \in ]1,+\infty[$}
We have the following error bounds and convergence result.
\begin{theo}
Suppose that $p \in ]1, +\infty[$. Let $u$ be a solution of~\eqref{neumann} with kernel $\K$ and data $(f,g)$. Let $u_n$ be a sequence of solutions to~\eqref{neumann} with kernels $\Kn$ and data $(f_n,g_n)$. Assume that $\K$ and $\Kn$ satisfy~\ref{assum:Kpos},~\ref{assum:Ksym} and $\K, \Kn \in L^{\infty,2}(\O^2)$, and that either one of the following holds:
\begin{enumerate}[label=(\alph*)]
\item $p \in ]1,2[$, $g, g_n \in L^2(\O)$, and $f, f_n \in L^1([0,T];L^2(\O))$; \label{assum:maina}
\item $p \geq 2$, $g, g_n \in L^{2(p-1)}(\O)$ and $f, f_n \in L^1([0,T];L^{2(p-1)}(\O))$; \label{assum:mainb}
\item $g, g_n \in L^{\infty}(\O)$ and $f, f_n \in L^1([0,T];L^{\infty}(\O))$. \label{assum:mainc}
\end{enumerate} 
Then, the following hold.
\begin{enumerate}[label=(\roman*)]
\item $u$ and $u_n$ are the unique solutions~\eqref{neumann} with respectively data $(f,g)$ and $(f_n,g_n)$. \label{item:maincontinuouswellposed}
\item We have the error estimate \label{item:maincontinuousbound}
\begin{multline}
\norm{u_n - u}_{C([0,T];L^2(\O))} \le 
\norm{g_n - g}_{L^2(\O)} + \norm{f_n - f}_{L^1([0,T];L^2(\O))} \\
+ CT
\begin{cases}
\norm{\Kn - \K}_{L^{\infty,2}(\O^2)}, 	& \text{under~\ref{assum:maina} or~\ref{assum:mainb}} \\
\norm{\Kn - \K}_{L^2(\O^2)}, 		& \text{under~\ref{assum:mainc}} 
\end{cases}
\label{eq:contbound}
\end{multline}
where $C$ is positive constant that may depend only on $p$, $g$ and $f$. 
\item Moreover, if~\ref{assum:mainc} holds, $\sup_{n \in \N} \abs{g_n(\bx)} < +\infty$ a.e. on $\O$ and $g_n \to g$ pointwise a.e. on $\O$, $\sup_{n \in \N} \abs{f_n(\bx,t))} < +\infty$ a.e. on $\O \times [0,T]$ and $f_n \to f$ pointwise a.e. on $\O \times [0,T]$, and the sequence $\seq{|\Kn|^2}$ is uniformly integrable over $\O^2$ and $\Kn \to \K$ pointwise a.e. on $\O^2$. Then 
\[
\lim_{n \to +\infty} \norm{u_n - u}_{C([0,T];L^2(\O))} = 0 .
\]
\label{item:maincontinuousconv}
\end{enumerate}
\label{thm:maincontinuous}
\end{theo}

\begin{rem}
Observe that since $L^{\infty}(\O) \subset L^2(\O)$ and $L^{2(p-1)}(\O) \subset L^2(\O)$ for $p \geq 2$, then the first two terms involved in~\eqref{eq:contbound} provide a non-trivial bound. Similarly, since $L^{\infty,2}(\O^2) \subset L^2(\O^2)$, the last term in the bound for case~\ref{assum:mainc} is also non-trivial.
%
In fact, both bounds in~\eqref{eq:contbound} can be summarized in one bound; the first one. However, the second bound for case~\ref{assum:mainc} is obviously sharper.
\end{rem}


\bpf{}
In the proof, $C$ is any positive constant that may depend solely on $p$ and $g$. 
\begin{enumerate}[label=(\roman*)]
\item Since $L^{\infty,2}(\O^2) \subset L^{\infty,1}(\O^2)$, assumption~\ref{assum:KL1} holds for both $\K$ and $\Kn$. We also have the embeddings
\begin{itemize}
\item $L^2(\O) \subset L^p(\O)$ under~\ref{assum:maina},
\item $L^{2(p-1)}(\O) \subset L^p(\O)$ under~\ref{assum:mainb}, and
\item $L^\infty(\O) \subset L^p(\O)$ under \ref{assum:mainc} . 
\end{itemize}
Thus $g,g_n \in L^p(\O)$ and $f,f_n \in L^1([0,T];L^p(\O))$. Existence and uniqueness of the solutions $u$ and $u_n$ in the sense of Definition~\ref{def:sol} is a consequence of Theorem~\ref{existence}.

\item Denote the error function $\xi_n(\bx,t) = u_n(\bx,t) - u(\bx,t)$, then from \eqref{neumann}, we have a.e.
\begin{equation}
\begin{split}
\frac{\partial \xi_n(\bx,t)}{\partial t }
&= - \pa{\plap{\Kn}{p}(u_n(\bx,t)) - \plap{\K}{p}(u(\bx,t))} + f_n(\bx,t) - f(\bx,t) \\
&= - \pa{\plap{\Kn}{p}(u_n(\bx,t)) - \plap{\Kn}{p}(u(\bx,t))} - \pa{\plap{\Kn}{p}(u(\bx,t)) - \plap{\K}{p}(u(\bx,t))}\\
&\quad + f_n(\bx,t) - f(\bx,t) .
\end{split}
\label{deriv}
\end{equation}
Multiplying both sides of~\eqref{deriv} by $\xi_n(\bx,t)$ and integrating, we get
\begin{equation}
\begin{split}
\frac{1}{2} \frac{\partial}{\partial t} \norm{\xi_n(\cdot,t)}_{L^2(\O)}^2 
 &= - \int_{\O} \pa{\plap{\Kn}{p} u_n(\bx,t) - \plap{\Kn}{p} u(\bx,t)}(u_n(\bx,t) - u(\bx,t)) d\bx \\
 & + \int_{\O^2} ( \Kn(\bx,\by) - \K(\bx,\by) ) \F(u(\by,t)-u(\bx,t)) \xi_n (\bx,t) d\bx d\by \\
 & + \int_{\O} \pa{f_n(\bx,t) - f(\bx,t)} \xi_n (\bx,t) d\bx .
\end{split}
\label{derivv}
\end{equation}
Since $g,g_n \in L^p(\O)$ and $f,f_n \in L^1([0,T];L^p(\O))$, $u_n(\cdot,t), u(\cdot,t) \in L^p(\O)$ for any $t \in [0,T]$ thanks to~\eqref{contraction}. We can then apply Proposition~\ref{prop:plap}\ref{itemmonotonep} with $h(x)=x$ to assert that the first term on the right-hand side of~\eqref{derivv} is nonpositive. Let us now bound the second term. 

\begin{itemize}
\item Case~\ref{assum:mainc}: in this case $\norm{u}_{C([0,T];L^\infty(\O))} \leq \norm{g}_{L^{\infty}(\O)} + \norm{f}_{L^1([0,T];L^{\infty}(\O))}$ thanks to~\eqref{contraction}, and we get from Cauchy-Schwartz inequality that
\begin{equation}
\begin{split}
& \abs { \int_{\O^2} (\Kn(\bx,\by) - \K(\bx,\by)) \F(u(\by,t) - u(\bx,t) ) \xi_n (\bx,t) d\bx d\by }\\
&\le 2^{p-1} \norm{u(\cdot,t)}_{L^{\infty}(\O)}^{p-1} \int_{\O^2} \abs{\Kn(\bx,\by) - \K(\bx,\by)} \abs{\xi_n(\bx,t)} d\bx d\by \\
&\le 2^{p-1} \pa{\norm{g}_{L^{\infty}(\O)}+\norm{f}_{L^1([0,T];L^{\infty}(\O))}}^{p-1} \norm{\Kn - \K}_{L^{2}(\O^2)} \norm{\xi_n(\cdot,t)}_{L^{2}(\O)} \\
&= C \norm{\Kn - \K}_{L^{2}(\O^2)} \norm{\xi_n(\cdot,t)}_{L^{2}(\O)} .
\end{split}
\label{estic}
\end{equation}

\item Case~\ref{assum:maina} or~\ref{assum:mainb}: applying again Cauchy-Schwartz inequality we obtain
\begin{equation*}
\hspace*{-1cm}
\begin{split}
& \abs {\int_{\O^2} (\Kn(\bx,\by) - \K(\bx,\by)) \F(u(\by,t) - u(\bx,t) ) \xi_n (\bx,t) d\bx d\by}\\
&\le \pa{\int_{\O^2} \abs{u(\by,t) - u(\bx,t)}^{2(p-1)} d\bx d\by}^{1/2} \pa{\int_{\O^2} |\Kn(\bx,\by) - \K(\bx,\by)|^2|\xi_n (\bx,t)|^2 d\bx d\by}^{1/2} \\
&= \pa{\int_{\O^2} \abs{u(\by,t) - u(\bx,t)}^{2(p-1)} d\bx d\by}^{1/2} \pa{\int_{\O} \pa{\int_{\O} |\Kn(\bx,\by) - \K(\bx,\by)|^2d\by} |\xi_n (\bx,t)|^2 d\bx}^{1/2} \\
&= \pa{\int_{\O^2} \abs{u(\by,t) - u(\bx,t)}^{2(p-1)} d\bx d\by}^{1/2} \norm{\Kn - \K}_{L^{\infty,2}(\O^2)} \norm{\xi_n(\cdot,t)}_{L^{2}(\O)} .
\end{split}
\end{equation*}
On the one hand, under~\ref{assum:maina}, Jensen's inequality applied to the concave function $x \in \R^+ \mapsto x^{p-1}$ entails
\begin{equation*}
\begin{split}
&\pa{\int_{\O^2} \abs{u(\by,t) - u(\bx,t)}^{2(p-1)} d\bx d\by}^{1/2} \\
&\leq \pa{\int_{\O^2} \abs{u(\by,t) - u(\bx,t)}^{2} d\bx d\by}^{(p-1)/2} \\
&\leq 2^{p-1}\norm{u(\cdot,t)}_{L^2(\O)}^{p-1} \leq 2^{p-1}\pa{\norm{g}_{L^2(\O)}+ \norm{f}_{L^1([0,T];L^2(\O))}}^{p-1} ,
\end{split}
\end{equation*} 
where we used~\eqref{contraction} in the last inequality. On the other hand, under~\ref{assum:mainb}, we have
\begin{equation*}
\begin{split}
\pa{\int_{\O^2} \abs{u(\by,t) - u(\bx,t)}^{2(p-1)} d\bx d\by}^{1/2} 
&\leq 2^{p-1}\norm{u(\cdot,t)}_{L^{2(p-1)}(\O)}^{p-1} \\
&\leq 2^{p-1}\pa{\norm{g}_{L^{2(p-1)}(\O)}+ \norm{f}_{L^1([0,T];L^{2(p-1)}(\O))}}^{p-1} .
\end{split}
\end{equation*} 
In turn, under either~\ref{assum:maina} or~\ref{assum:mainb}, we have the bound
\begin{equation}
\begin{split}
& \abs {\int_{\O^2} (\Kn(\bx,\by) - \K(\bx,\by)) \F(u(\by,t) - u(\bx,t) ) \xi_n (\bx,t) d\bx d\by}\\
&\le C \norm{\Kn - \K}_{L^{\infty,2}(\O^2)} \norm{\xi_n(\cdot,t)}_{L^{2}(\O)} .
\end{split}
\label{estiab}
\end{equation}
\end{itemize}
Inserting~\eqref{estic} and~\eqref{estiab} into~\eqref{derivv}, ignoring the first term which is non-positive as argued above, and using Cauchy-Schwartz inequality on the last term, we obtain
\begin{equation*}
\frac{\partial}{\partial t} \norm{\xi_n(\cdot,t)}_{L^2(\O)} 
\le 
\norm{f_n(\cdot,t) - f(\cdot,t)}_{L^2(\O)} + 
\begin{cases}
C \norm{\Kn - \K}_{L^{\infty,2}(\O^2)}, & \text{under \ref{assum:maina} or \ref{assum:mainb}} \\
C \norm{\Kn - \K}_{L^2(\O^2)}, & \text{under \ref{assum:mainc}}.
\end{cases}
\end{equation*}
Integrating this inequality on $[0,t]$ and taking the supremum over $t \in [0,T]$, we get \eqref{eq:contbound}.

\item By assumptions on $\seq{\Kn}$, we are in position to apply the Vitali convergence theorem~\cite[p.~133]{Rudin} in $L^2(\O^2)$ to get that $\norm{\Kn - \K}_{L^2(\O^2)} \to 0$ as $n \to +\infty$. We have by assumption that the sequence $\seq{g_n}$ is dominated by a constant function. The latter is obviously in $L^2(\O)$ since $|\O| < +\infty$. It then follows from the dominated convergence theorem that $\norm{g_n - g}_{L^2(\O)} \to 0$ as $n \to +\infty$. We now turn to the sequence $f_n$. We have
\[
\norm{f_n - f}_{L^1([0,T];L^2(\O))} \leq T^{1/2} \norm{f_n - f}_{L^2([0,T];L^2(\O))} = T^{1/2}\norm{f_n - f}_{L^2(\O \times [0,T])} .
\] 
Arguing as for $g_n$, using our assumptions, entails again that $\norm{f_n - f}_{L^1([0,T];L^2(\O))} \to 0$ as $n \to +\infty$. Passing to the limit in the second inequality of~\eqref{eq:contbound}, the claim follows.
\end{enumerate}
\epf

In the case where the kernel takes the form $\K(\bx,\by)=\J(\bx-\by)$, we have the following consequence of Theorem~\ref{thm:maincontinuous}. 
\begin{cor}
Suppose that $p \in ]1, +\infty[$. Let $u$ be a solution of \eqref{neumann} with kernel $\K(\bx,\by)=\J(\bx-\by)$ and data $(f,g)$. Let $u_n$ be a sequence of solutions to \eqref{neumann} with kernels $\Kn(\bx,\by)=\Jn(\bx-\by)$ and data $(f_n,g_n)$. Assume that $\J$ and $\Jn$ satisfy~\ref{assum:Jpos},~\ref{assum:Jsym} and $\J, \Jn \in L^2(\O-\O)$, and that either one of~\ref{assum:maina},~\ref{assum:mainb} or~\ref{assum:mainc} in Theorem~\ref{thm:maincontinuous} holds.
Then, the following hold.
\begin{enumerate}[label=(\roman*)]
\item $u$ and $u_n$ are the unique solutions of the corresponding evolution problems. 
\item We have the error estimate
\begin{equation}
\norm{u_n - u}_{C([0,T];L^2(\O))} \le \norm{g_n - g}_{L^2(\O)} + \norm{f_n - f}_{L^1([0,T];L^2(\O))} + CT\norm{\Jn - \J}_{L^2(\O - \O)} ,
\label{eq:contboundJ}
\end{equation}
where $C$ is positive constant that may depend only on $p$, $g$ and $f$. 
\item Moreover, if the sequence $\seq{|\Jn|^2}$ is uniformly integrable over $\O-\O$, $\Jn \to \J$ pointwise a.e. on $\O-\O$, $g_n \to g$ pointwise a.e. on $\O$, $f_n \to f$ pointwise a.e. on $\O \times [0,T]$, and either one of the following holds:
\begin{enumerate}[label=(\alph*')]
\item $p \in ]1,2[$, $\seq{|g_n|^2}$ (resp. $\seq{|f_n|^2}$) is uniformly integrable over $\O$ (resp. $\O \times [0,T]$); \label{assum:mainJap}
\item $p \geq 2$, $\seq{|g_n|^{2(p-1)}}$ (resp. $\seq{|f_n|^{2(p-1)}}$) is uniformly integrable over $\O$ (resp. $\O \times [0,T]$); \label{assum:mainJbp}
\item $\sup_{n \in \N} \abs{g_n(\bx)} < +\infty$ a.e. on $\O$ and $\sup_{n \in \N} \abs{f_n(\bx,t)} < +\infty$ a.e. on $\O \times [0,T]$. \label{assum:mainJcp}
\end{enumerate} 
Then 
\[
\lim_{n \to +\infty} \norm{u_n - u}_{C([0,T];L^2(\O))} = 0 .
\]
\end{enumerate}
\label{cor:maincontinuousJ}
\end{cor}

\bpf{}
\begin{enumerate}[label=(\roman*)]
\item We argue in the same way as in the proof Theorem~\ref{thm:maincontinuous} since $L^2(\O-\O) \subset L^1(\O-\O)$ implies that assumption~\ref{assum:JL1} holds for both $\J$ and $\Jn$.
\item The error bound~\eqref{eq:contboundJ} is a specialization of~\eqref{eq:contbound} since
\[
\int_{\O} |\Kn(\bx,\by)-\K(\bx,\by)|^2 d\by = \int_{\O-\bx} |\Jn(\bz)-\J(\bz)|^2 d\bz \leq \norm{\Jn-\J}_{L^2(\O-\O)}^2 .
\]
Thus
\[
\norm{\Kn-\K}_{L^2(\O^2)} \leq \norm{\Kn-\K}_{L^{\infty,2}(\O^2)} \leq \norm{\Jn-\J}_{L^2(\O-\O)}.
\]
\item Case~\ref{assum:mainJap} follows from the Vitali convergence theorem applied to $\Jn$, $g_n$ and $f_n$. The latter argument also applies to case \ref{assum:mainJbp} since $L^{2(p-1)}(\O-\O) \subset L^2(\O-\O)$, $L^{2(p-1)}(\O) \subset L^2(\O)$ and $L^{2(p-1)}(\O \times [0,T]) \subset L^1([0,T];L^2(\O))$. Case~\ref{assum:mainJcp} uses the Vitali convergence theorem on $\Jn$ and the dominated convergence theorem on $g_n$ and $f_n$ as argued in the proof of Theorem~\ref{thm:maincontinuous}\ref{item:maincontinuousconv}.
\end{enumerate}
\epf{}

\begin{rem}
At this stage, we only relied on the monotonicity property of $\plap{K}{p}$ in Proposition~\ref{prop:plap}\ref{itemmonotonep} to get our bounds. One may then wonder if the stronger notions of monotonicity established in Proposition~\ref{prop:plap}\ref{itemstrongmonotonep} can yield bounds better than~\eqref{eq:contboundJ}. We answer this question positively by (slightly) improving the dependence on $T$ for $p \in ]1,2]$ but at the price of more stringent assumptions on $\J$. For this, we embark from~\eqref{derivv}, bound all terms as in the proof of Theorem~\ref{thm:maincontinuous}, use Proposition~\ref{prop:plap}\ref{itemstrongmonotonep} and that $L^2(\O) \subset L^p(\O)$ in this case to get
\begin{multline*}
\frac{1}{2} \frac{\partial}{\partial t} \norm{\xi_n(\cdot,t)}_{L^2(\O)}^2 
 + C_1 \int_{\O^2} \J(\bx-\by) \abs{\gradNL\xi_n(\bx,\by)}^2 d\by d\bx
 \leq \\
 \pa{C \norm{\Jn - \J}_{L^2(\O - \O)} + \norm{f_n(\cdot,t) - f(\cdot,t)}_{L^2(\O)}}\norm{\xi_n(\cdot,t)}_{L^{2}(\O)} ,
\end{multline*}
for two positive constants $C,C_1$ (in the following $C_i$ is a positive constant). Assume in addition that $\J$ is compactly supported and $\J(0) > 0$. One can then invoke the Poincar\'e inequality~\cite[Proposition~4.1]{mazonjmaa} to show that
\[
C_2 \norm{\xi_n(\cdot,t) - \int_{\O}\xi_n(\bx,t) d\bx}_{L^2(\O)}^2 \leq \int_{\O^2} \J(\bx-\by) \abs{\gradNL\xi_n(\bx,\by)}^2 d\by d\bx .
\]
Thus
\[
\frac{1}{2}\norm{\xi_n(\cdot,t)}_{L^2(\O)}^2 \leq \norm{\xi_n(\cdot,t) - \int_{\O}\xi_n(\bx,t) d\bx}_{L^2(\O)}^2 + \pa{\int_{\O}\xi_n(\bx,t) d\bx}^2 .
\]
Altogether, we arrive at
\begin{multline*}
\frac{1}{2} \frac{\partial}{\partial t} \norm{\xi_n(\cdot,t)}_{L^2(\O)}^2 
 + \frac{C_1C_2}{2} \norm{\xi_n(\cdot,t)}_{L^2(\O)}^2
 \leq \\
 \pa{C \norm{\Jn - \J}_{L^2(\O - \O)} + \norm{f_n(\cdot,t) - f(\cdot,t)}_{L^2(\O)}}\norm{\xi_n(\cdot,t)}_{L^{2}(\O)} + C_1C_2 \pa{\int_{\O}\xi_n(\bx,t) d\bx}^2 .
\end{multline*}
By integrating~\eqref{neumann}, it is easy to see by applying Proposition~\ref{prop:plap}\ref{itemstrongmonotonep} and~\ref{itemmonotonep} with $h(x)=1$ that the solution of~\eqref{neumann} preserves the total mass in $\O$, whence we deduce
\[
\int_{\O}\xi_n(\bx,t) d\bx = \int_{\O} (g_n(\bx)-g(\bx)) + \int_{0}^t\int_{\O}(f_n(\bx,s)-f(\bx,s)) d\bx ds .
\]
If $(f,g)$ and $(f_n,g_n)$ have the same mass, we get
\begin{multline*}
\frac{1}{2} \frac{\partial}{\partial t} \norm{\xi_n(\cdot,t)}_{L^{2}(\O)}^2 
 + \frac{C_1C_2}{2} \norm{\xi_n(\cdot,t)}_{L^2(\O)}^2
 \leq \\
 \pa{C\norm{\Jn - \J}_{L^2(\O - \O)} + \norm{f_n(\cdot,t) - f(\cdot,t)}_{L^2(\O)}}\norm{\xi_n(\cdot,t)}_{L^{2}(\O)},
\end{multline*}
and therefore
\begin{align*}
\frac{\partial}{\partial t} \norm{\xi_n(\cdot,t)}_{L^{2}(\O)} + \frac{C_1C_2}{2} \norm{\xi_n(\cdot,t)}_{L^2(\O)}
\leq \pa{C\norm{\Jn - \J}_{L^2(\O - \O)} + \norm{f_n(\cdot,t) - f(\cdot,t)}_{L^2(\O)}} .
\end{align*}
Applying Gronwall's lemma yields the estimate
\begin{multline*}
\norm{u_n(\bx,t)-u(\bx,t)}_{L^{2}(\O)} \leq \norm{f_n - f}_{L^1([0,T];L^2(\O))} + \exp(-C_1C_2 t/2) \norm{g_n-g}_{L^2(\O)} \\
+ \frac{2C}{C_1C_2}(1-\exp(-C_1C_2 t/2)) \norm{\Jn - \J}_{L^2(\O - \O)} .
\end{multline*}
This bound is clearly better than~\eqref{eq:contboundJ}. In turn,
\begin{multline*}
\norm{u_n-u}_{C([0,T];L^{2}(\O))} \leq \norm{f_n - f}_{L^1([0,T];L^2(\O))} +
\max\pa{\norm{g_n-g}_{L^2(\O)},\frac{2C}{C_1C_2}\norm{\Jn - \J}_{L^2(\O - \O)}} . 
\end{multline*}
The same reasoning as above can be carried out to sharpen the error bounds for the discrete problems in Section~\ref{sec:discrete}. Nevertheless, this will not be detailed further in this work.
\end{rem}

\subsection{The case $p=1$}
We now turn to the case $p=1$.
\begin{theo}
Let $u$ be a solution of~\eqref{neumann} for $p=1$ with kernel $\K$ and data $(f,g)$. Let $u_n$ be a sequence of solutions to~\eqref{neumann} for $p=1$ with kernels $\Kn$ and data $(f_n,g_n)$. Assume that $\K$ and $\Kn$ satisfy~\ref{assum:Kpos} and~\ref{assum:Ksym}, that $\K, \K_n \in L^{\infty,2}(\O^2)$, $g, g_n \in L^2(\O)$ and $f, f_n \in L^1([0,T];L^2(\O))$.
Then, the following hold.
\begin{enumerate}[label=(\roman*)]
\item $u$ and $u_n$ are the unique solutions in the sense of Definition~\ref{def:solp1} of the corresponding evolution problems. \label{item:maincontinuouswellposedp1}
\item We have the error estimate \label{item:maincontinuousboundp1}
\begin{equation}
\norm{u_n - u}_{C([0,T];L^2(\O))} \le \norm{g_n - g}_{L^2(\O)} + \norm{f_n - f}_{L^1([0,T];L^2(\O))} + T\norm{\Kn - \K}_{L^2(\O^2)} .
\label{eq:contboundp1}
\end{equation}
\item Moreover, if $\Kn \to \K$ pointwise a.e. on $\O^2$, $g_n \to g$ pointwise a.e. on $\O$, $f_n \to f$ pointwise a.e. on $\O \times [0,T]$, and $\seq{|\Kn|^2}$ is uniformly integrable over $\O^2$, $\seq{|g_n|^2}$ is uniformly integrable on $\O$, and $\seq{|f_n|^2}$ is uniformly integrable on $\O \times [0,T]$. Then 
\[
\lim_{n \to +\infty} \norm{u_n - u}_{C([0,T];L^2(\O))} = 0 .
\]
\label{item:maincontinuousconvp1}
\end{enumerate}
\label{thm:maincontinuousp1}
\end{theo}

 
\bpf{}
\begin{enumerate}[label=(\roman*)]
\item Existence and uniqueness of $u$ and $u_n$ follow from Theorem~\ref{existencep1} where we argue as in Theorem~\ref{thm:maincontinuous}\ref{item:maincontinuouswellposed} since $g, g_n \in L^2(\O) \subset L^1(\O)$ and $\K, \K_n \in L^{\infty,2}(\O^2) \subset L^{\infty,1}(\O^2)$. 
\item Denote the error function $\xi_n(\bx,t) = u_n(\bx,t) - u(\bx,t)$, then from Definition~\ref{def:solp1}, we have a.e.
\begin{equation}
\begin{split}
\frac{\partial \xi_n(\bx,t)}{\partial t } 
&= \int_{\O} \pa{\Kn(\bx,\by) w_n(\bx,\by,t) d\by - \K(\bx,\by)w(\bx,\by,t)} d\by + f_n(\bx,t) - f(\bx,t) \\
&= \int_{\O} \Kn(\bx,\by) \pa{w_n(\bx,\by,t) - w(\bx,\by,t)} d\by + 
\int_{\O} (\Kn(\bx,\by)-\K(\bx,\by)) w(\bx,\by,t) d\by \\
& + f_n(\bx,t) - f(\bx,t) ,
\end{split}
\label{derivp1}
\end{equation}
where $w$ (resp. $w_n$) is the subgradient function associated to $u$ (resp. $u_n$) as in Definition~\ref{def:solp1}. Multiplying both sides of~\eqref{derivp1} by $\xi_n(\bx,t)$ and integrating, we get
\begin{equation}
\begin{split}
\frac{1}{2} \frac{\partial}{\partial t} \norm{\xi_n(\cdot,t)}_{L^2(\O)}^2 
&= \int_{\O^2} \Kn(\bx,\by) \pa{w_n(\bx,\by,t) - w(\bx,\by,t)} \xi_n(\bx,t) d\bx d\by \\
&+ \int_{\O^2} (\Kn(\bx,\by)-\K(\bx,\by)) w(\bx,\by,t)\xi_n(\bx,t) d\bx d\by \\
&+ \int_{\O} \pa{f_n(\bx,t) - f(\bx,t)} \xi_n (\bx,t) d\bx .
\end{split}
\label{derivvp1}
\end{equation}
In view of the monotonicity claim in Proposition~\ref{prop:1lap}\ref{itemmonotonep1}, we have
\[
\int_{\O^2} \Kn(\bx,\by) \pa{w_n(\bx,\by,t) - w(\bx,\by,t)} \xi_n(\bx,t) d\bx d\by \leq 0 .
\]
Let us turn to bounding the second term. We have by the Cauchy-Schwartz inequality and that $\norm{w}_{L^\infty(\O^2 \times ]0,T[)} \leq 1$,
\begin{equation}
\begin{split}
& \abs {\int_{\O^2} (\Kn(\bx,\by) - \K(\bx,\by)) w(\bx,\by,t) \xi_n (\bx,t) d\bx d\by }\\
&\le \int_{\O^2} \abs{\Kn(\bx,\by) - \K(\bx,\by)} \abs{\xi_n(\bx,t)} d\bx d\by \\
&\le \norm{\Kn - \K}_{L^{2}(\O^2)} \norm{\xi_n(\cdot,t)}_{L^{2}(\O)} .
\end{split}
\label{estip1}
\end{equation}

Inserting~\eqref{estip1} into~\eqref{derivvp1}, ignoring the first term which is non-positive as argued above, and using Cauchy-Schwartz inequality on the last term, we obtain
\begin{equation*}
\frac{\partial}{\partial t} \norm{\xi_n(\cdot,t)}_{L^2(\O)} 
\leq \norm{f_n(\cdot,t) - f(\cdot,t)}_{L^2(\O)} + \norm{\Kn - \K}_{L^{2}(\O^2)} .
\end{equation*}
Integrating this inequality on $[0,t]$ and taking the supremum over $t \in [0,T]$, we get~\eqref{eq:contboundp1}.

\item We argue again using the Vitali convergence theorem since $\K, \Kn \in L^{\infty,2}(\O^2) \subset L^2(\O^2)$ and $L^1([0,T];L^2(\O)) \subset L^2(\O \times [0,T])$.
\end{enumerate}
\epf

The following corollary is immediate in the same vein as Corollary~\ref{cor:maincontinuousJ}.
\begin{cor}
Let $u$ be a solution of~\eqref{neumann} for $p=1$ with kernel $\K(\bx,\by)=\J(\bx-\by)$ and data $(f,g)$. Let $u_n$ be a sequence of solutions to~\eqref{neumann} for $p=1$ with kernels $\Kn(\bx,\by)=\Jn(\bx-\by)$ and data $(f_n,g_n)$. Assume that $\J$ and $\Jn$ satisfy~\ref{assum:Jpos},~\ref{assum:Jsym} and $\J, \Jn \in L^2(\O-\O)$, that $g, g_n \in L^2(\O)$ and $f, f_n \in L^1([0,T];L^2(\O))$. Then, the following hold.
\begin{enumerate}[label=(\roman*)]
\item $u$ and $u_n$ are the unique solutions in the sense of Definition~\ref{def:solp1} of the corresponding evolution problems. 
\item We have the error estimate
\begin{equation}
\norm{u_n - u}_{C([0,T];L^2(\O))} \le \norm{g_n - g}_{L^2(\O)} + \norm{f_n - f}_{L^1([0,T];L^2(\O))} + T\norm{\Jn - \J}_{L^2(\O - \O)} .
\label{eq:contboundJp1}
\end{equation}
\item Moreover, if $\Jn \to \J$ pointwise a.e. on $\O - \O$, $g_n \to g$ pointwise a.e. on $\O$, $f_n \to f$ pointwise a.e. on $\O \times [0,T]$, and $\seq{|\Jn|^2}$ is uniformly integrable over $\O-\O$, $\seq{|g_n|^2}$ is uniformly integrable on $\O$, and $\seq{|f_n|^2}$ is uniformly integrable on $\O \times [0,T]$. Then 
\[
\lim_{n \to +\infty} \norm{u_n - u}_{C([0,T];L^2(\O))} = 0 .
\]
\end{enumerate}
\label{cor:maincontinuousJp1}
\end{cor}

\section{Error bounds for the discrete problem}
\label{sec:discrete}

Let $\bK \in \R^{n^d \times n^d}$ and $\bg \in \R^{n^d}$ be discrete approximations of, respectively, the kernel $\K$ and initial data $g$ in~\eqref{neumann}, on a regular mesh of size $\deltan$. Typically, one can take $\bK = \projn K$ and $\bg = \projn g$. 
For $1<p<\infty$, the discrete $p$-Laplacian operator with kernel $\bK$ is
\begin{align*}
\plapdisc{\bK}{p}: \bu \in \R^{n^d} \mapsto - \sum_{{\bj} \in [n]^d} h_{\bj} \bK_{\bi\bj} \abs{\bu_{\bj} - \bu_{\bi}}^{p-2} (\bu_{\bj} - \bu_{\bi}) = -\sum\limits_{{\bj} \in [n]^d} h_{\bj} \bK_{\bi\bj} \Psi(\bu_{\bj} - \bu_{\bi}) .
\end{align*}
In the same way, we define the discrete $1$-Laplacian operator as the set-valued operator $\plapdisc{\bK}{1}: \R^{n^d} \to 2^{\R^{n^d}}$ such that ${\bs \eta} \in \plapdisc{\bK}{1} \bu$ if and only if
\[
{\bs \eta}_\bi = - \sum_{{\bj} \in [n]^d} h_{\bj} \bK_{\bi\bj} \bw_{\bi\bj} , 
\]
where $\norm{\bw}_{\infty} \leq 1$, $\bw_{\bi\bj}=-\bw_{\bj\bi}$, and 
\[
\bw_{\bi\bj} \in \partial |\cdot| (\bu_{\bj}-\bu_{\bi}) .
\]

By construction, we have the following simple lemma whose proof is immediate.
\begin{lem}
For any $\bK \in \R^{n^d \times n^d}$ and $\bu \in \R^{n^d}$, the following holds:
\begin{enumerate}[label=(\roman*)]
\item If $1 < p < +\infty$, \label{item:injplap}
\[
\injn \plapdisc{\bK}{p}(\bu) = \plap{\injn \bK}{p}(\injn\bu).
\]
\item If $p=1$, \label{item:injplap1}
\[
\injn {\bs \eta}(\bx) = -\int_{\O} \injn\K(\bx,\by) \injn\bw(\bx,\by) d\by, \qwhereq \injn\bw(\bx,\by) \in \partial |\cdot| (\injn\bu(\by)-\injn\bu(\bx)) .
\]
Moreover, $\norm{\injn\bw}_{L^{\infty}(\O^2)} \leq 1$ and $\injn\bw(\bx,\by)=-\injn\bw(\by,\bx)$.
\end{enumerate}
\label{lem:injplap}
\end{lem}

\subsection{The semi-discrete problem}
\label{semidiscrete}
\paragraph{Case $p \in ]1,+\infty[$:} 
We start with the case $1 < p < +\infty$ and consider the space semi-discretization of~\eqref{neumann},
\begin{equation}\tag{\textrm{${\mathcal{P}}^{\text{SD}}_{p}$}}
\begin{cases}
\frac{\partial }{\partial t} \bu(t) = -\plapdisc{\bK}{p} \bu(t) + \bff(t), & t > 0, \\
\bu(0) = \bg .
\end{cases}
\label{neumannsemidisc}
\end{equation}
where $\bu: t \in \R^+ \mapsto \bu(t) \in \R^{n^d}$ and similarly for $\bff$.

Our aim is to compare the solutions of problems~\eqref{neumann} and~\eqref{neumannsemidisc}. The solution of~\eqref{neumannsemidisc} being discrete in space, we consider its continuum space extensions of $\bu$ and $\bff$ on $\O$ for any $t > 0$ as
\begin{equation}\label{eq:spaceextuf}
u_n(\bx,t) = (\injn \bu(t))(\bx) \qandq f_n(\bx,t) = (\injn \bff(t))(\bx) .
\end{equation}

\begin{theo}
Suppose that $p \in ]1, +\infty[$. Let $u$ be a solution of~\eqref{neumann} with kernel $\K$ and data $(f,g)$, and $\bu$ that of~\eqref{neumannsemidisc} with $\bK=\projn \K$, $\bg=\projn g$ and $\bff(t)=\projn f(\cdot,t)$ for $t \in [0,T]$. Let $u_n$ and $f_n$ as defined in~\eqref{eq:spaceextuf}. Assume that $\K$ satisfies~\ref{assum:Kpos},~\ref{assum:Ksym} and $\K \in L^{\infty,2}(\O^2)$, and that $g$ and $f$ satisfy either one of the conditions~\ref{assum:maina},~\ref{assum:mainb} or~\ref{assum:mainc} in Theorem~\ref{thm:maincontinuous}. Then, the following hold.
\begin{enumerate}[label=(\roman*)]
\item $u$ and $u_n$ are the unique solutions of~\eqref{neumann} with data respectively $(f,g)$ and $(f_n,\injn\projn g)$. 
\item We have the error estimate
\begin{multline}
\norm{u_n - u}_{C([0,T];L^2(\O))} \le 
\norm{\injn\projn g - g}_{L^2(\O)} + \norm{f_n - f}_{L^1([0,T];L^2(\O))} \\
+ CT
\begin{cases}
\norm{\injn\projn \K - \K}_{L^{\infty,2}(\O^2)}, 	& \text{under~\ref{assum:maina}-\ref{assum:mainb}} \\
\norm{\injn\projn \K - \K}_{L^2(\O^2)}, 		& \text{under~\ref{assum:mainc}} 
\end{cases}
\label{eq:semidiscbound}
\end{multline}
where $C$ is positive constant that depends only on $p$, $g$ and $f$. 
\item If, moreover, $g \in L^{\infty}(\O) \cap \Lip(s,L^2(\O))$, $\K \in \Lip(s,L^2(\O^2))$ and $f(\cdot,t) \in L^\infty(\O) \cap \Lip(s,L^2(\O))$ for every $t \in [0,T]$, then
\begin{equation}
\norm{u_n - u}_{C([0,T];L^2(\O))} \leq C(1+T)\deltan^{s} ,
\label{eq:semidiscrate}
\end{equation}
where $C$ is positive constant that depends only on $p$, $g$, $f$, $\K$, $s$.
\end{enumerate}
\label{thm:mainsemidiscrete}
\end{theo}

\bpf{}
\begin{enumerate}[label=(\roman*)]
\item Existence and uniqueness of $u$ were proved in Theorem~\ref{thm:maincontinuous}\ref{item:maincontinuouswellposed}. We also see that $\injn\bK$ verifies~\ref{assum:Kpos} and~\ref{assum:Ksym}. Using Lemma~\ref{lem:projinj}, we have $\injn \bg \in L^p(\O)$, $f_n \in L^1([0;T],L^p(\O))$ and $\injn\bK \in L^{\infty,2}(\O^2) \subset L^{\infty,1}(\O^2)$, and thus $\injn\bK$ fulfills~\ref{assum:KL1}. In view of Lemma~\ref{lem:injplap}\ref{item:injplap}, it follows from~\eqref{neumannsemidisc} that the function $u_n$ satisfies~\eqref{neumann} with kernel $\injn\bK$ and data $(f_n,\injn\bg)$. Existence and uniqueness of $u_n$ then follow from Theorem~\ref{existence}.

\item The claim is a specialization of~\eqref{eq:contbound} in Theorem~\ref{thm:maincontinuous}\ref{item:maincontinuousbound}.

\item As $\K \in L^{\infty,2}(\O^2) \subset L^2(\O^2) $, we insert the estimate \eqref{eq:lipspaceapprox} (see Lemma~\ref{lem:spaceapprox}) in the second bound of \eqref{eq:semidiscbound}.
\end{enumerate}
\epf{}

\paragraph{Case $p=1$:} 
We now turn to the case $p=1$, and consider the evolution problem
\begin{equation}\tag{\textrm{${\mathcal{P}}^{\text{SD}}_{1}$}}
\begin{cases}
\frac{\partial }{\partial t} \bu(t) = -{\bs\eta}(t) + \bff(t), & t > 0, \\
\bu(0) = \bg ,
\end{cases}
\label{neumannsemidiscp1}
\end{equation}
where 
\[
{\bs\eta}_\bi(t) = -\sum_{{\bj} \in [n]^d} h_{\bj} \bK_{\bi\bj} \sign(\bu_{\bj} - \bu_{\bi}), \quad \text{and thus} \quad {\bs\eta}(t) \in \plapdisc{\bK}{1} \bu(t) .
\]

\begin{theo}
Let $u$ be a solution of~\eqref{neumann} for $p=1$ with kernel $\K$ and data $(f,g)$, and $\bu$ is that of~\eqref{neumannsemidiscp1} with $\bK=\projn \K$, $\bg=\projn g$ and $\bff(t)=\projn f(\cdot,t)$ for $t \in [0,T]$. Let $u_n$ and $f_n$ as defined in~\eqref{eq:spaceextuf}. Assume that $\K$ satisfies~\ref{assum:Kpos},~\ref{assum:Ksym} and $\K \in L^{\infty,2}(\O^2)$, and that $g \in L^2(\O)$ and $f \in L^1([0,T];L^2(\O))$. Then, the following hold.
\begin{enumerate}[label=(\roman*)]
\item $u$ and $u_n$ are the unique solutions in the sense of Definition~\ref{def:solp1} of the corresponding evolution problems.
\item We have the error estimate
\begin{equation}
\norm{u_n - u}_{C([0,T];L^2(\O))} \le \norm{\injn\projn g - g}_{L^2(\O)} + \norm{f_n - f}_{L^1([0,T];L^2(\O))} + T\norm{\injn\projn K - \K}_{L^2(\O^2)} .
\label{eq:semidiscboundp1}
\end{equation}
\item If, moreover, $g \in \Lip(s,L^2(\O))$, $\K \in \Lip(s,L^2(\O^2))$ and $f(\cdot,t) \in \Lip(s,L^2(\O))$ for every $t \in [0,T]$, then
\begin{equation}
\norm{u_n - u}_{C([0,T];L^2(\O))} \leq C(1+T)\deltan^{s} ,
\label{eq:semidiscratep1}
\end{equation}
where $C$ is positive constant that depends only on $p$, $g$, $f$, $\K$ and $s$.
\end{enumerate}
\label{thm:mainsemidiscretep1}
\end{theo}

\bpf{}
\begin{enumerate}[label=(\roman*)]
\item Existence and uniqueness of $u$ were proved in Theorem~\ref{thm:maincontinuousp1}\ref{item:maincontinuouswellposedp1}. In addition, $\injn\bK$ verifies~\ref{assum:Kpos} and~\ref{assum:Ksym}. Using Lemma~\ref{lem:projinj}, $\injn \bg \in L^2(\O) \subset L^1(\O)$, $f_n \in L^1([0,T];L^2(\O)) \subset L^1([0,T];L^1(\O))$ and $\injn\bK \in L^{\infty,2}(\O^2) \subset L^{\infty,1}(\O^2)$, and thus $\injn\bK$ fulfills~\ref{assum:KL1}. By virtue of Lemma~\ref{lem:injplap}\ref{item:injplap1}, $u_n$, the space continuum extension of $\bu$, will satisfy~\eqref{neumann} with kernel $\injn\bK$ and data $(f_n,\injn\bg)$. Existence and uniqueness of $u_n$ in the sense of Definition~\ref{def:solp1} follow from Theorem~\ref{existencep1}.

\item This claim is a specialization of~\eqref{eq:contboundp1} in Theorem~\ref{thm:maincontinuousp1}\ref{item:maincontinuousboundp1}.

\item Insert the estimate~\eqref{eq:lipspaceapprox} in \eqref{eq:semidiscboundp1}.
\end{enumerate}
\epf{}

\subsection{The totally discrete problem}
\label{totallydiscrete}

We establish in this section error bounds for fully discrete (in time and space) approximations of~\eqref{neumann}. For that, let $0 < t_1 < t_2 < \cdots < t_{N-1} < t_N=T$ be a partition (not necessarily equispaced) of $[0,T]$. Let $\tauhm \eqdef \abs{t_k - t_{k-1}}$ and denote $\tau = \max\limits_{k \in [N]} \tauh$.

\subsubsection{Forward/Explicit Euler discretization}\label{subsec:forwardeuler}
 
\paragraph{Case $p \in ]1,2]$:}
We start with $p \in ]1,2]$ and consider a totally discrete problem with forward/explicit Euler scheme in time,
\begin{equation}\tag{\textrm{${\mathcal{P}}^{\text{TDF}}_{p}$}}
\begin{cases}
\displaystyle{\frac{\bu^{k} - \bu^{k-1}}{\tauhm}} = -\plapdisc{\bK}{p} \bu^{k-1} + \bff, & k \in [N], \\
\bu^0 = \bg ,
\end{cases}
\label{neumannfulldiscforward}
\end{equation}
where $\bu^k, \bff \in \R^{n^d}$. We have implicitly assumed that $\bff$ does not depend on time, which is a standard assumption in the context of explicit discretization.

Since our aim is to compare the solutions of problems~\eqref{neumann} and~\eqref{neumannfulldiscforward}, we introduce the following continuum extensions in space and/or time of $\acc{\bu^k}_{k \in [N]}$ as
\begin{align*}
u_n^k &= \injn \bu^k, k \in [N], \qandq f_n = \injn \bff, \\
\unlin(\bx,t) &= \frac{t_{k} - t}{\tauhm} u_n^{k-1}(\bx) + \frac{t-t_{k-1}}{\tauhm} u_n^k(\bx) , \quad (\bx,t) \in \O \times ]t_{k-1},t_k], k \in [N], \\
\uncst(\bx,t) &= \sum_{k=1}^N u_n^{k-1}(\bx) \chi_{]t_{k-1},t_k]} (t), \quad (\bx,t) \in \O \times ]0,T].
\end{align*} 

Then, in the same vein as Lemma~\ref{lem:injplap}, it is easy to see that~\eqref{neumannfulldiscforward} is equivalent to the following evolution problem 
\begin{equation}
\begin{cases}
\frac{\partial}{\partial t} \unlin(\bx,t)=-\plap{\injn\bK}{p}\uncst(\bx,t) + f_n(\bx), & (\bx,t) \in \O \times ]0,T], \\
\unlin(\bx,0) = \injn\bg(\bx), \quad \bx \in \O .
\end{cases}
\label{neumannfulldiscforwardinterp}
\end{equation}

Before turning to the consistency result, we collect some useful estimates.
\begin{lem} 
Consider problem~\eqref{neumannfulldiscforward} with kernel $\bK$, data $(\bff,\bg)$ and variable step-size $\tauh \leq 2C \norm{\plap{\injn\bK}{p}u_n^k - f_n}_{L^2(\O)}^{\frac{2-p}{p-1}}$, where $C$ is the constant in Proposition~\ref{prop:plap}\ref{itemholderp}. Assume that $\injn\bg \in L^2(\O)$ and $\injn \bK$ satisfies~\ref{assum:Kpos},~\ref{assum:Ksym} and~\ref{assum:KL1}. Suppose also that for each $n \in \N$, $\bff$ is such that~\eqref{neumannfulldiscforward} has a stationary solution $\bu^\star$ and that $\sup_{n \in \N}\norm{\injn\bg - \injn\bu^\star}_{L^2(\O)} < +\infty$. Then
\[
\text{$\uncst(\cdot,t) \in L^2(\O)$}, \forall t \in [0,T], \qandq \sup_{t \in [0,T], n \in \N} \norm{\uncst(\cdot,t) - \injn\bu^\star}_{L^2(\O)} < +\infty .
\]
\label{lem:estimatesuncstforward}
\end{lem}

\begin{rem}
{~}\vspace*{0.1cm}
\begin{enumerate}[label=(\arabic*)]
\item Condition on the time-step $\tauh$ can be seen as an abstract non-linear CFL condition. It is better than the one in \cite{Hafienne_FE_Nonlocal_plaplacian_evolu} since we here exploit the H\"older continuity of $\plap{\injn\bK}{p}$ on $L^2(\O)$ for $p \in ]1,2]$, see Proposition~\ref{prop:plap}\ref{itemholderp}. For $p=2$, where $\plap{\injn\bK}{2}$ is linear Lipschitz continuous operator on $L^2(\O)$, the condition reads $\tauh \leq 2C$. Such condition for explicit time-discretization of evolution problems with accretive and Lipschitz-continuous operators is known, see e.g.,~\cite{Savare06}. It is also consistent with known convergence results for finding zeros of co-called co-coercive operators on Hilbert spaces~\cite{bauschke2011convex}. \label{item1}
\item The assumption on $\bff$ and $\bK$ imply that $f_n \in L^2(\O)$. Indeed,~\eqref{eq:plapholder} entails
\[
\norm{f_n}_{L^2(\O)} = \norm{\plap{\injn\bK}{p}(\injn\bu^\star)}_{L^2(\O)} \leq 2^p C_2\norm{\K}_{L^{\infty,1}(\O^2)}^2\norm{\injn\bu^\star}_{L^2(\O)}^{p-1} .
\]
\item The assumption made on $\bff$ is trivially true when $\bff = \bf{0}$ since $\bf{0}$ is a stationary solution int this case. In turn, using Lemma~\ref{lem:projinj}, one can see that the uniform boundedness conditions on $\bg$ and $\bK$ are fulfilled if $\bg = \projn g$ and $\bK = \projn\K$, where $g \in L^2(\O)$ and $\K$ satisfies~\ref{assum:Kpos}-\ref{assum:KL1}. \label{item2}
\end{enumerate}
\label{rem:estimatesuncst}
\end{rem}

\bpf{}
We show the claim by an induction argument. Since $\plap{\injn\bK}{p}(\injn\bu^\star) = f_n$, we have
\begin{align*}
&\norm{u_n^{1} - \injn\bu^\star}_{L^2(\O)}^2 \\
&=\norm{\injn \bg - \injn \bu^\star}_{L^2(\O)}^2 - 2\tau_0\int_{\O} \pa{\plap{\injn\bK}{p}(\injn\bg)(\bx)-f_n(\bx)}\pa{\injn\bg(\bx) - \injn \bu^\star} d\bx\\
&\qquad  + \tau_0^2 \norm{\plap{\injn\bK}{p}(\injn\bg) - f_n}_{L^2(\O)}^2 \\
&=\norm{\injn \bg - \injn \bu^\star}_{L^2(\O)}^2 - 2\tau_0\int_{\O} \pa{\plap{\injn\bK}{p}(\injn\bg)(\bx)-\plap{\injn\bK}{p}(\injn\bu^\star)(\bx)}\pa{\injn\bg(\bx) - \injn\bu^\star} d\bx \\
&\qquad + \tau_0^2 \norm{\plap{\injn\bK}{p}(\injn\bg) - f_n}_{L^2(\O)}^2 .
\end{align*}
By assumption on $\bg$, $\bu^\star$ and $\tauh$, we can invoke Proposition~\ref{prop:plap}\ref{itemholderp} to get
\begin{align*}
&\norm{u_n^{1} - \injn\bu^\star}_{L^2(\O)}^2 \\
&\leq \norm{\injn \bg - \injn\bu^\star}_{L^2(\O)}^2 - 2C\tau_0 \norm{\plap{\injn\bK}{p}(\injn\bg)-f_n}_{L^2(\O)}^{p/(p-1)} + \tau_0^2 \norm{\plap{\injn\bK}{p}(\injn\bg)-f_n}_{L^2(\O)}^2 \\
&\leq \norm{\injn \bg - \injn \bu^\star}_{L^2(\O)}^2 + \tau_0\norm{\plap{\injn\bK}{p}(\injn\bg) - f_n}_{L^2(\O)}^2\pa{2C\norm{\plap{\injn\bK}{p}(\injn\bg) - f_n}_{L^2(\O)}^{(2-p)/(p-1)} - \tau_0} \\
&\leq \norm{\injn \bg - \injn \bu^\star}_{L^2(\O)}^2 .
\end{align*}
Suppose now that, for any $k > 1$,
\[
\norm{u_n^{k} - \injn \bu^\star}_{L^2(\O)}^2 \leq \norm{\injn \bg - \injn \bu^\star}_{L^2(\O)}^2 ,
\]
and thus $u_n^{k} \in L^2(\O)$. We can then use Proposition~\ref{prop:plap}\ref{itemholderp} as above to see that
\begin{align*}
&\norm{u_n^{k+1} - \injn \bu^\star}_{L^2(\O)}^2 \\
&\leq \norm{\injn \bg - \injn \bu^\star}_{L^2(\O)}^2 - \tauh\norm{\plap{\injn\bK}{p}(u_n^k) - f_n}_{L^2(\O)}^2\pa{2C\norm{\plap{\injn\bK}{p}(u_n^k) - f_n}_{L^2(\O)}^{(2-p)/(p-1)} - \tauh} \\
&\leq \norm{\injn \bg - \injn \bu^\star}_{L^2(\O)}^2 .
\end{align*}
Thus the sequence $\acc{\norm{u_n^{k}}_{L^2(\O)}}_{k \in \N}$ is bounded, and so is $\norm{\uncst(\cdot,t)}_{L^2(\O)}$ for $t \in [0,T]$ by its definition. We also have
\begin{align*}
\sup_{t \in [0,T], n \in \N} \norm{\uncst(\cdot,t) - \injn\bu^\star}_{L^2(\O)} 
&= \sup_{(n,N) \in \N^2, k \in [N]} \norm{u_n^{k} - \injn\bu^\star}_{L^2(\O)} \leq \sup_{n \in \N} \norm{\injn \bg - \injn \bu^\star}_{L^2(\O)} < +\infty .
\end{align*}

\epf{}

\begin{lem} 
In addition to the assumptions of Lemma~\ref{lem:estimatesuncstforward}, suppose that $\sup_{n \in \N}\norm{\injn\bK}_{L^{\infty,1}(\O^2)} < +\infty$. Then
\begin{equation*}
\sup_{t \in [0,T], n \in \N}\norm{\unlin(\cdot,t)- \uncst(\cdot,t)}_{L^2(\O)} \leq C \tau,
\end{equation*}
where $C$ is a positive constant that does not depend on $(n,N,T)$.
\label{lem:unlincstforward}
\end{lem}
\bpf{}
It is easy to see that for $t \in ]t_{k-1},t_k]$, $k \in \N$,
\begin{align*}
\norm{\unlin(\cdot,t) - \uncst(\cdot,t)}_{L^2(\O)} 
&= (t - t_{k-1}) \norm{\frac{u_n^k -u_n^{k-1}}{\tauhm}}_{L^2(\O)} \\
&= (t - t_{k-1}) \norm{\plap{\injn \bK}{p}u_n^{k-1} - f_n}_{L^2(\O)} \\
&= (t - t_{k-1}) \norm{\plap{\injn \bK}{p}u_n^{k-1} - \plap{\injn \bK}{p}\injn\bu^\star}_{L^2(\O)} \\
&\leq \tau \norm{\plap{\injn \bK}{p}u_n^{k-1} - \plap{\injn \bK}{p}\injn\bu^\star}_{L^2(\O)} = \tau \norm{\plap{\injn \bK}{p}\uncst(\cdot,t) - \plap{\injn \bK}{p}\injn\bu^\star}_{L^2(\O)} .
\end{align*}
As $\plap{\injn\bK}{p}$ is H\"older continuous on $L^2(\O)$ with exponent $p-1$, see~\eqref{eq:plapholder}, we get
\begin{align*}
\norm{\unlin(\cdot,t) - \uncst(\cdot,t)}_{L^2(\O)} 
&\leq \tau 2^{p/2} C_2^{1/2}\norm{\K}_{L^{\infty,1}(\O^2)}\norm{\uncst(\cdot,t) - \injn\bu^\star}_{L^2(\O)}^{p-1} .
\end{align*}
We then take the supremum over $t$ and $n$, and use Lemma~\ref{lem:estimatesuncstforward} to conclude.
\epf


We are now in position to state the error bound for the totally discrete problem~\eqref{neumannfulldiscforward}.

\begin{theo}
Suppose that $p \in ]1, 2]$. Let $u$ be a solution of~\eqref{neumann} with kernel $\K$ and data $(f,g)$ where $f$ is time-independent, and $\acc{\bu^k}_{k \in [N]}$ is the sequence generated by~\eqref{neumannfulldiscforward} with $\bK=\projn \K$, $\bg=\projn g$, $\bff=\projn f$ and $\tau_k$ as prescribed in Lemma~\ref{lem:estimatesuncstforward}. Assume that $\K$ satisfies~\ref{assum:Kpos},~\ref{assum:Ksym} and $\K \in L^{\infty,2}(\O^2)$, and that $f,g$ belong either to $L^2(\O)$ or $L^\infty(\O)$. Then, the following hold.
\begin{enumerate}[label=(\roman*)]
\item $u$ is the unique solution of~\eqref{neumann}, $\acc{\bu^k}_{k \in [N]}$ is uniquely defined and $\acc{\norm{\injn\bu^k}_{L^2(\O)}}_{k \in [N]}$ is bounded (uniformly in $n$ when $\bff=\bf{0}$).
\item We have the error estimate
\begin{multline}
\sup_{k \in [N], t \in ]t_{k-1},t_k]}\norm{\injn\bu^{k-1} - u(\cdot,t)}_{L^2(\O)} 
\leq \exp\pa{T/2}\Bigg(\norm{\injn\projn g - g}_{L^2(\O)} \\
+ CT^{1/2} \pa{\tau^{1/(3-p)} + \norm{f_n - f}_{L^2(\O)} +
\begin{cases}
\norm{\injn\projn \K - \K}_{L^{\infty,2}(\O^2)} & g \in L^2(\O) \\
\norm{\injn\projn \K - \K}_{L^2(\O^2)},	& g \in L^\infty(\O)
\end{cases}}\Bigg) .
\label{eq:fulldiscboundforwardcst}
\end{multline}
for $\tau$ sufficiently small, where $C$ is positive constant that depends only on $p$, $g$, $f$ and $\K$.
\item If, moreover, $f,g \in L^{\infty}(\O) \cap \Lip(s,L^2(\O))$ and $\K \in \Lip(s,L^2(\O^2))$, then
\begin{equation}
\sup_{k \in [N], t \in ]t_{k-1},t_k]}\norm{\injn\bu^{k-1} - u(\cdot,t)}_{L^2(\O)} \leq C\exp(T/2)\pa{(1+T^{1/2})\deltan^{s} + T^{1/2}\tau^{1/(3-p)}},
\label{eq:fulldiscreteforwardratecst}
\end{equation}
for $\tau$ sufficiently small, where $C$ is positive constant that depends only on $p$, $g$, $f$, $\K$ and $s$.
\end{enumerate}
\label{thm:mainfulldiscreteforward}
\end{theo}

\bpf{}
In the proof, $C$ is any positive constant that may depend only on $p$, $g$, $f$, $\K$ and/or $s$, and that may be different at each line.
\begin{enumerate}[label=(\roman*)]
\item Existence and uniqueness of $u$ were proved in Theorem~\ref{thm:maincontinuous}\ref{item:maincontinuouswellposed}. The claimed well-posedness of the sequence $\acc{\bu^k}_{k \in [N]}$ is a consequence of Lemma~\ref{lem:estimatesuncstforward} and Remark~\ref{rem:estimatesuncst}\ref{item2}.

\item Denote $\xilin(\bx,t) = \unlin(\bx,t) - u(\bx,t)$, $\xicst(\bx,t) = \uncst(\bx,t) - u(x,t)$, $g_n = \injn\projn g$ and $\Kn = \injn\projn\K$. We thus have a.e.
\begin{equation*}
\begin{split}
\frac{\partial \xilin(\bx,t)}{\partial t } 
&= - \pa{\plap{\Kn}{p}(\uncst(\bx,t)) - \plap{\K}{p}(u(\bx,t))} + (f_n(\bx)-f(\bx)) \\
&= - \pa{\plap{\Kn}{p}(\uncst(\bx,t)) - \plap{\Kn}{p}(u(\bx,t))} - \pa{\plap{\Kn}{p}(u(\bx,t)) - \plap{\K}{p}(u(\bx,t))} + (f_n(\bx)-f(\bx)) .
\end{split}
\end{equation*}
Multiplying both sides by $\xilin(\bx,t)$, integrating and rearranging the terms, we get
\begin{equation}
\begin{split}
\frac{1}{2} \frac{\partial}{\partial t} \norm{\xilin(\cdot,t)}_{L^2(\O)}^2 
 &= - \int_{\O} \pa{\plap{\Kn}{p} \uncst(\bx,t) - \plap{\Kn}{p} u(\bx,t)}(\uncst(\bx,t) - u(\bx,t)) d\bx \\
 &- \int_{\O} \pa{\plap{\Kn}{p} u(\bx,t) - \plap{\K}{p} u(\bx,t)}\xilin(\bx,t) d\bx \\
 & - \int_{\O} \pa{\plap{\Kn}{p} \uncst(\bx,t) - \plap{\Kn}{p} u(\bx,t)}\pa{\unlin(\bx,t) - \uncst(\bx,t)} d\bx \\
 & + \int_{\O} \pa{f_n(\bx) - f(\bx)} \xilin (\bx,t) d\bx .
\end{split}
\label{eq:derivfullforward}
\end{equation}
Since $f,g \in L^p(\O)$ in both cases, so is $u(\cdot,t)$ thanks to \eqref{contraction}. We also have $\uncst(\cdot,t) \in L^2(\O) \subset L^p(\O)$ by Lemma~\ref{lem:estimatesuncstforward}. We are then in position to use Proposition~\ref{prop:plap}\ref{itemmonotonep} with $h(x)=x$ to assert that the first term on the right-hand side of~\eqref{eq:derivfullforward} is nonpositive. Let us now bound the second term. 

Similarly to the estimates~\eqref{estiab} and~\eqref{estic} in the proof of Theorem~\ref{thm:maincontinuous}, and using Young inequality, we have
\begin{align*}
&\abs{\int_{\O} \pa{\plap{\Kn}{p} u(\bx,t) - \plap{\K}{p} u(\bx,t)}\xilin(\bx,t) d\bx} \\
&\leq
\begin{cases}
C \norm{\Kn - \K}_{L^{\infty,2}(\O^2)}\norm{\xilin(\cdot,t)}_{L^2(\O)}, & g \in L^2(\O) \\
C \norm{\Kn - \K}_{L^2(\O^2)}\norm{\xilin(\cdot,t)}_{L^2(\O)}, 		& g \in L^\infty(\O),
\end{cases} \\
&\leq \frac{1}{6}\norm{\xilin(\cdot,t)}_{L^2(\O)}^2 +
\begin{cases}
C \norm{\Kn - \K}_{L^{\infty,2}(\O^2)}^2, 	& g \in L^2(\O) \\
C \norm{\Kn - \K}_{L^2(\O^2)}^2, 		& g \in L^\infty(\O).
\end{cases}
\end{align*}
For the third term in~\eqref{eq:derivfullforward}, we invoke Lemma~\ref{lem:unlincstforward} to get
\begin{align*}
& \abs{\int_{\O} \pa{\plap{\Kn}{p} \uncst(\bx,t) - \plap{\Kn}{p} u(\bx,t)}\pa{\unlin(\bx,t) - \uncst(\bx,t)} d\bx} \\
&\leq \norm{\plap{\Kn}{p} \uncst(\cdot,t)-\plap{\Kn}{p}u(\cdot,t)}_{L^2(\O)}\norm{\unlin(\cdot,t) - \uncst(\cdot,t)}_{L^2(\O)} \\
&\leq C\norm{\plap{\Kn}{p} \uncst(\cdot,t)-\plap{\Kn}{p}u(\cdot,t)}_{L^2(\O)}\tau .
\end{align*}
We then use the fact that $\plap{\injn\bK}{p}$ is H\"older continuous on $L^2(\O)$ with exponent $p-1$, see~\eqref{eq:plapholder}, to obtain
\begin{align*}
\norm{\plap{\Kn}{p} \uncst(\cdot,t)-\plap{\Kn}{p}u(\cdot,t)}_{L^2(\O)}
\leq C\norm{\xicst(\cdot,t)}_{L^2(\O)}^{p-1} 
\leq C\pa{\norm{\xilin(\cdot,t)}_{L^2(\O)}^{p-1} + \tau^{p-1}} ,
\end{align*}
where we used Lemma~\ref{lem:unlincstforward} to go from $\xicst$ to $\xilin$, and that $p \in ]1,2]$. It then follows by Cauchy-Schwartz inequality that
\begin{align*}
& \abs{\int_{\O} \pa{\plap{\Kn}{p} \uncst(\bx,t) - \plap{\Kn}{p} u(\bx,t)}\pa{\unlin(\bx,t) - \uncst(\bx,t)} d\bx} \\
&\leq C\pa{\norm{\xilin(\cdot,t)}_{L^2(\O)}^{p-1}\tau + \tau^p} \\
&\leq \frac{1}{6}\norm{\xilin(\cdot,t)}_{L^2(\O)}^2 + C(\tau^{2/(3-p)} + \tau^p) .
\end{align*}
Using Young inequality to bound the last term in~\eqref{eq:derivfullforward}, and combining the bounds on the three other terms, we have shown that
\begin{multline*}
\frac{\partial}{\partial t} \norm{\xilin(\cdot,t)}_{L^2(\O)}^2 
\leq \norm{\xilin(\cdot,t)}_{L^2(\O)}^2 + C\Bigg(\tau^{2/(3-p)} + \tau^p + \norm{f_n - f}_{L^2(\O)}^2 \\
+
\begin{cases}
\norm{\Kn - \K}_{L^{\infty,2}(\O^2)}^2, & g \in L^2(\O) \\
\norm{\Kn - \K}_{L^2(\O^2)}^2,	& g \in L^\infty(\O)
\end{cases}\Bigg) .
\end{multline*}

Using the Gronwall's lemma and taking the square-root, we get
\begin{multline}
\norm{\unlin - u}_{C([0,T];L^2(\O))}
\leq \exp\pa{T/2}\Bigg(\norm{\injn\projn g - g}_{L^2(\O)} \\
+ CT^{1/2} \pa{\tau^{1/(3-p)} + \tau^{p/2} + \norm{f_n - f}_{L^2(\O)} +
\begin{cases}
\norm{\injn\projn \K - \K}_{L^{\infty,2}(\O^2)}	& g \in L^2(\O) \\
\norm{\injn\projn \K - \K}_{L^2(\O^2)},		& g \in L^\infty(\O)
\end{cases}}\Bigg) .
\label{eq:fulldiscboundforward}
\end{multline}
Since $1/2 < 1/(3-p) \leq p/2$ for $p \in ]1,2]$ the dependence on $\tau$ scales as $O(\tau^{1/(3-p)})$ for $\tau$ sufficiently small (or $N$ large enough). Inserting \eqref{eq:fulldiscboundforward} into
\begin{align}
\sup_{k \in [N], t \in ]t_{k-1},t_k]}\norm{u_n^{k-1} - u(\cdot,t)}_{L^2(\O)} 
= \norm{\uncst - u}_{C([0,T];L^2(\O))} 
\leq \norm{\unlin - u}_{C([0,T];L^2(\O))} + C\tau ,
\label{eq:bnduncstulin}
\end{align}
completes the proof of the error bound.

\item Plug~\eqref{eq:lipspaceapprox} into~\eqref{eq:fulldiscboundforwardcst}.
\end{enumerate}
\epf

\begin{rem}
Error bounds in $L^p(\O)$ were derived in~\cite{Hafienne_FE_Nonlocal_plaplacian_evolu} for forward Euler discretization. Their rate is better than ours and is provided for the range $p \in ]1,+\infty[$. Unfortunately, we believe that their proof contains invalid arguments that can be fixed but only for $p \in ]1,2]$.
\end{rem}

\paragraph{Case $p=1$:}
We now turn to the case $p=1$, and consider the discrete system
\begin{equation}\tag{\textrm{${\mathcal{P}}^{\text{TDF}}_{1}$}}
\begin{cases}
\displaystyle{\frac{\bu^{k} - \bu^{k-1}}{\tauhm}} = -{\bs\eta}^{k-1} + \bff, & k \in [N], \\
\bu^0 = \bg .
\end{cases}
\label{neumannfulldiscforwardp1}
\end{equation}
where
\[
{\bs\eta}^{k} = -\sum_{{\bj} \in [n]^d} h_{\bj} \bK_{\bi\bj} \sign(\bu_{\bj}^{k} - \bu_{\bi}^{k}), \quad \text{and thus} \quad {\bs\eta}^{k} \in \plapdisc{\bK}{1} \bu^{k} .
\]
We consider the continuum extensions in space and/or time of $\acc{\bu^k}_{k \in [N]}$ as before, namely $u_n^k$, $\unlin$ and $\uncst$, $f_n = \injn \bff$, and the space-time continuum extension of $\acc{{\bs\eta}^k}_{k \in [N]}$
\[
\etancst(\bx,t) = \sum_{k=1}^N (\injn {\bs\eta}^{k-1})(\bx) \chi_{]t_{k-1},t_k]} (t) = -\int_{\O} \injn\bK(\bx,\by) \sign(\uncst(\by,t)-\uncst(\bx,t)), \quad (\bx,t) \in \O \times ]0,T] .
\]
In view of Lemma~\ref{lem:injplap}, these extensions satisfy the evolution problem 
\begin{equation}
\begin{cases}
\frac{\partial}{\partial t} \unlin(\bx,t) = -\etancst(\bx,t) + f_n(\bx), & (\bx,t) \in \O \times ]0,T], \\
\unlin(\bx,0) = \injn\bg(\bx), \quad \bx \in \O ,
\end{cases}
\label{neumannfulldiscforwardinterpp1}
\end{equation}
and
\[
\etancst(\bx,t) \in \plap{\injn\bK}{1}\uncst(\bx,t) .
\]

We have the following counterpart estimates of Lemma~\ref{lem:estimatesuncstforward}.
\begin{lem} 
Consider problem~\eqref{neumannfulldiscforwardp1} with kernel $\bK$, data $(\bff,\bg)$ and variable step-size $$\tauh = \frac{\alpha_k}{\max\pa{\norm{\injn{\bs\eta}^{k} - f_n}_{L^2(\O)},1}},\qwhereq \sum_{k \in \N}\alpha_{k}^2 < + \infty.$$ Assume that $\injn\bg \in L^2(\O)$ and $\injn \bK$ satisfies~\ref{assum:Kpos}-\ref{assum:Ksym} and~\ref{assum:KL1}. Suppose also that for each $n \in \N$, $\bff$ is such that~\eqref{neumannfulldiscforward} has a stationary solution $\bu^\star$ and that $\sup_{n \in \N}\norm{\injn\bg - \injn\bu^\star}_{L^2(\O)} < +\infty$. Then
\[
\text{$\uncst(\cdot,t) \in L^2(\O)$}, \forall t \in [0,T], \qandq \sup_{t \in [0,T], n \in \N} \norm{\uncst(\cdot,t) - \injn\bu^\star}_{L^2(\O)} < +\infty .
\]
\label{lem:estimatesuncstforwardp1}
\end{lem}

\begin{rem}
The condition on the time-step $\tauh$ is reminescent of subgradient descent and has been used in~\cite{Hafienne_FE_Nonlocal_plaplacian_evolu}. The assumptions on $(\bff,\bg,\bK)$ are again verified when $\bff = 0$, $\bg = \projn g$ and $\bK = \projn\K$, where $g \in L^2(\O)$ and $\K$ satisfies~\ref{assum:Kpos}-\ref{assum:KL1}.
\end{rem}

\bpf{}
Define the series $s_k \eqdef \sum_{i=0}^k \alpha_i^2$. As in Lemma~\ref{lem:estimatesuncstforward}, we proceed by induction using the monotonicity of the $1$-Laplacian (Proposition~\ref{prop:1lap}\ref{itemmonotonep1}). Indeed, since $f_n \in \plap{\injn\bK}{p}(\injn\bu^\star)$, we have
\begin{multline*}
\norm{u_n^{1} - \injn\bu^\star}_{L^2(\O)}^2 = \norm{\injn \bg - \injn \bu^\star}_{L^2(\O)}^2 \\
- 2\tau_0\int_{\O} \pa{\plap{\injn\bK}{p}(\injn\bg)(\bx)-\plap{\injn\bK}{p}(\injn\bu^\star)(\bx)}\pa{\injn\bg(\bx) - \injn\bu^\star} d\bx + \alpha_0^2 .
\end{multline*}
By assumption on $\bg$, $\bu^\star$, we can invoke Proposition~\ref{prop:1lap}\ref{itemmonotonep1} to get
\begin{align*}
\norm{u_n^{1} - \injn\bu^\star}_{L^2(\O)}^2 \leq \norm{\injn \bg - \injn\bu^\star}_{L^2(\O)}^2 + s_0.
\end{align*}
Suppose now that, for any $k > 1$,
\[
\norm{u_n^{k} - \injn \bu^\star}_{L^2(\O)}^2 \leq \norm{\injn\bg - \injn \bu^\star}_{L^2(\O)}^2 + s_{k-1} ,
\]
and thus $u_n^{k} \in L^2(\O)$. We can then invoke again Proposition~\ref{prop:1lap}\ref{itemmonotonep1} to see that
\begin{align*}
&\norm{u_n^{k+1} - \injn\bu^\star}_{L^2(\O)}^2 \\
&= \norm{u_n^k - \injn \bu^\star}_{L^2(\O)}^2 - 2\tau_k\int_{\O} \pa{\plap{\injn\bK}{p}(u_n^k)(\bx)-\plap{\injn\bK}{p}(\injn\bu^\star)(\bx)}\pa{u_n^k(\bx) - \injn\bu^\star} d\bx + \alpha_k^2 \\
&\leq \norm{\injn\bg - \injn \bu^\star}_{L^2(\O)}^2 + s_k .
\end{align*}
This shows that for all $k \in \N$,
\begin{align*}
\norm{u_n^k - \injn\bu^\star}_{L^2(\O)}^2 \leq \norm{\injn\bg - \injn \bu^\star}_{L^2(\O)}^2 + s_{\infty} ,
\end{align*}
and thus $\acc{\norm{\injn\bu^k}_{L^2(\O)}}_{k \in [N]}$ is bounded. In turn, so is $\norm{\uncst(\cdot,t)}_{L^2(\O)}$ for $t \in [0,T]$ by its definition. Moreover,
\begin{align*}
\sup_{t \in [0,T], n \in \N} \norm{\uncst(\cdot,t) - \injn\bu^\star}_{L^2(\O)} 
&= \sup_{(n,N) \in \N^2, k \in [N]} \norm{u_n^{k} - \injn\bu^\star}_{L^2(\O)} \\
&\leq \sup_{n \in \N}\norm{\injn \bg - \injn\bu^\star}_{L^2(\O)} + s_{\infty}^{1/2} < +\infty .
\end{align*}
\epf{}

We also have the following analogue of Lemma~\ref{lem:unlincstforward}.
\begin{lem} 
In addition to the assumptions of Lemma~\ref{lem:estimatesuncstforwardp1}, suppose that $\sup_{n \in \N}\norm{\injn\bK}_{L^{\infty,1}(\O^2)} < +\infty$. Then
\begin{equation*}
\sup_{t \in [0,T], n \in \N}\norm{\unlin(\cdot,t)- \uncst(\cdot,t)}_{L^2(\O)} \leq C \tau,
\end{equation*}
where $C$ is a positive constant that does not depend on $(n,N,T)$.
\label{lem:unlincstforwardp1}
\end{lem}

\bpf{}
Arguing as the beginning of Lemma~\ref{lem:unlincstforward}, we get for any $t \in ]t_{k-1},t_k]$, $k \in \N$,
\begin{align*}
\norm{\unlin(\cdot,t) - \uncst(\cdot,t)}_{L^2(\O)} 
&\leq \tau \norm{\etancst(\bx,t) - f_n}_{L^2(\O)} .
\end{align*}
By H\"older inequality, we have
\begin{align*}
\norm{\etancst(\bx,t)}_{L^2(\O)}^2
&=\int_{\O}\aabs{\int_{\O}\injn\bK(\bx,\by)\sign(\uncst(\by,t)-\uncst(\bx,t))d\by}^2d\bx \\
&\leq\int_{\O}\pa{\int_{\O}\injn\bK(\bx,\by)d\by}^2d\bx 
\leq \norm{\injn\bK}_{L^{\infty,1}(\O^2)}^2 .
\end{align*}
Thee same bound also holds on $\norm{f_n}_{L^2(\O)}$. We then take the supremum over $t$ and $n$ to conclude.
\epf

\begin{theo}
Let $u$ be a solution of~\eqref{neumann} with kernel $\K$ and data $(f,g)$ where $f$ is time-independent, and $\acc{\bu^k}_{k \in [N]}$ is the sequence generated by~\eqref{neumannfulldiscforwardp1} with $\bK=\projn \K$, $\bg=\projn g$, $\bff=\projn f$ and $\tau_k$ as prescribed in Lemma~\ref{lem:estimatesuncstforwardp1}. Assume that $\K$ satisfies~\ref{assum:Kpos},~\ref{assum:Ksym} and $\K \in L^{\infty,2}(\O^2)$, and that $f, g \in L^2(\O)$. Then, the following hold.
\begin{enumerate}[label=(\roman*)]
\item $u$ is the unique solution of \eqref{neumann}, $\acc{\bu^k}_{k \in [N]}$ is uniquely defined and $\acc{\norm{\injn\bu^k}_{L^2(\O)}}_{k \in [N]}$ is bounded (uniformly in $n$ when $\bff=\bf{0}$).
\item We have the error estimate
\begin{multline}
\sup_{k \in [N], t \in ]t_{k-1},t_k]}\norm{\injn\bu^{k-1} - u(\cdot,t)}_{L^2(\O)}
\leq \exp\pa{T/2}\Bigg(\norm{\injn\projn g - g}_{L^2(\O)} \\
+ CT^{1/2} \pa{\tau^{1/2} + \norm{f_n - f}_{L^2(\O)} + \norm{\injn\projn \K - \K}_{L^2(\O^2)}}\Bigg)
\label{eq:fulldiscboundforwardp1}
\end{multline}
where $C$ is positive constant that depends only on $\K$.
\item If, moreover, $f, g \in \Lip(s,L^2(\O))$ and $\K \in \Lip(s,L^2(\O^2))$, then
\begin{equation}
\sup_{k \in [N], t \in ]t_{k-1},t_k]}\norm{\injn\bu^{k-1} - u(\cdot,t)}_{L^2(\O)} \leq C\exp(T/2)\pa{(1+T^{1/2})\deltan^{s} + T^{1/2}\tau^{1/2}},
\label{eq:fulldiscrateforward}
\end{equation}
where $C$ is positive constant that depends only on $g$, $f$, $\K$ and $s$.
\end{enumerate}
\label{thm:mainfulldiscreteforwardp1}
\end{theo}

\bpf{}
$C$ is any positive constant that may depend only on $g$, $f$, $\K$ and $s$, and that may be different at each line. We use the same notation as in the proof of Theorem~\ref{thm:mainfulldiscreteforward}.
\begin{enumerate}[label=(\roman*)]
\item Existence and uniqueness of $u$ were proved in Theorem~\ref{thm:maincontinuous}\ref{item:maincontinuouswellposed}. Well-posedness of $\acc{\bu^k}_{k \in [N]}$ follows from Lemma~\ref{lem:estimatesuncstforwardp1} and Remark~\ref{rem:estimatesuncst}\ref{item2}.

\item We have
\begin{multline*}
\frac{\partial \xilin(\bx,t)}{\partial t} 
=\int_{\O} \Kn(\bx,\by) \pa{\wncst(\bx,\by,t) - w(\bx,\by,t)} d\by \\
+ \int_{\O} (\Kn(\bx,\by)-\K(\bx,\by)) w(\bx,\by,t) d\by + (f_n(\bx) - f(\bx)) ,
\end{multline*}
where $w$ is the subgradient function associated to $u$ (see Definition~\ref{def:solp1}), and $\wncst(\bx,\by,t) = \sign(\uncst(\by,t)-\uncst(\bx,t))$. Multiplying both sides by $\xilin(\bx,t)$, integrating and rearranging the terms, we get
\begin{equation}
\begin{split}
\frac{1}{2} \frac{\partial}{\partial t} \norm{\xilin(\cdot,t)}_{L^2(\O)}^2 
 &= \int_{\O^2} \Kn(\bx,\by) \pa{\wncst(\bx,\by,t) - w(\bx,\by,t)}(\uncst(\bx,t) - u(\bx,t)) d\bx d\by \\
 & + \int_{\O^2} (\Kn(\bx,\by)-\K(\bx,\by)) w(\bx,\by,t)\xilin(\bx,t) d\bx d\by \\
 & + \int_{\O} \Kn(\bx,\by) \pa{\wncst(\bx,\by,t) - w(\bx,\by,t)}\pa{\unlin(\bx,t) - \uncst(\bx,t)} d\bx d\by \\
 & + \int_{\O} \pa{f_n(\bx) - f(\bx)} \xilin (\bx,t) d\bx .
\end{split}
\label{eq:derivfullforwardp1}
\end{equation}
As $u(\cdot,t) \in L^1$ and $\uncst(\cdot,t) \in L^2(\O) \subset L^1(\O)$ by Lemma~\ref{lem:estimatesuncstforwardp1}, the monotonicity claim in Proposition~\ref{prop:1lap}\ref{itemmonotonep1} yields that the first term in~\eqref{eq:derivfullforwardp1} is nonpositive. The second and third terms can be easily bounded as
\begin{align*}
\abs{\int_{\O^2} (\Kn(\bx,\by)-\K(\bx,\by)) w(\bx,\by,t)\xilin(\bx,t) d\bx d\by}
&\leq \norm{\Kn-\K}_{L^2(\O^2)}\norm{\xilin(\cdot,t)}_{L^2(\O)} \\
&\leq \frac{1}{4}\norm{\xilin(\cdot,t)}_{L^2(\O)}^2 + \norm{\Kn-\K}_{L^2(\O^2)}^2 .
\end{align*}
and the third term using Lemma~\ref{lem:unlincstforwardp1}
\[
\abs{\int_{\O} \Kn(\bx,\by) \pa{\wncst(\bx,\by,t) - w(\bx,\by,t)}\pa{\unlin(\bx,t) - \uncst(\bx,t)} d\bx d\by} 
\leq 2 \norm{\K}_{L^{\infty,2}(\O^2)}^2 \tau .
\]
Bounding the last term by Young inequality, we obtain
\begin{equation*}
\frac{\partial}{\partial t} \norm{\xilin(\cdot,t)}_{L^2(\O)}^2 
\leq \norm{\xilin(\cdot,t)}_{L^2(\O)}^2 + 2\norm{f_n-f}_{L^2(\O)}^2 + 2\norm{\Kn-\K}_{L^2(\O^2)}^2 + C\tau .
\end{equation*}
Using the Gronwall's lemma and~\eqref{eq:bnduncstulin}, we get the claimed bound.

\item Insert~\eqref{eq:lipspaceapprox} into~\eqref{eq:fulldiscboundforwardp1}.
\end{enumerate}
\epf{}

\subsubsection{Backward/Implicit Euler discretization}\label{subsec:backwardeuler}
Forward Euler discretization was able to deal only with $p \in [1,2]$. For backward Euler discretization, we will tackle $p \in ]1,+\infty[$.

We consider the fully discrete problem with backward Euler time scheme
\begin{equation}\tag{\textrm{${\mathcal{P}}^{\text{TDB}}_{p}$}}
\begin{cases}
\displaystyle{\frac{\bu^{k} - \bu^{k-1}}{\tauhm}} = -\plapdisc{\bK}{p} \bu^{k} + \bff^k, & k \in [N], \\
\bu^0 = \bg ,
\end{cases}
\label{neumannfulldiscbackward}
\end{equation}
where $\bu^k, \bff^k \in \R^{n^d}$. This can also be written equivalently as
\begin{equation*}
\bu^{k} = J_{\tauhm\plapdisc{\bK}{p}}(\bu^{k-1} + \tauhm \bff^k) .
\end{equation*}
This is known as the proximal iteration, and is at the heart of so-called mild solutions as well as existence and uniqueness of solutions to~\eqref{neumann} through the nonlinear semigroups theory~\cite{CrandallLiggett71,Benilan72,kobayashi,BenilanCrandall91}. Denoting as before $u_n^k = \injn \bu^k$ and $f_n^k = \injn \bff^k$ the space continuum extensions of $\bu^k$ and $\bff^k$, we also have
\begin{equation*}
u_n^{k} = J_{\tauhm\plap{\injn\bK}{p}}(u_n^{k-1} + \tauhm f_n^k) .
\end{equation*}
We also let the time-space continuum extensions
\begin{align*}
\unlin(\bx,t) &= \frac{t_{k} - t}{\tauhm} u_n^{k-1}(\bx) + \frac{t-t_{k-1}}{\tauhm} u_n^k(\bx) , \quad (\bx,t) \in \O \times ]t_{k-1},t_k], k \in [N], \\
\uncst(\bx,t) &= \sum_{k=1}^N u_n^{k}(\bx) \chi_{]t_{k-1},t_k]} (t) \qandq \fncst(\bx,t) = \sum_{k=1}^N f_n^{k}(\bx) \chi_{]t_{k-1},t_k]} (t), \quad (\bx,t) \in \O \times ]0,T].
\end{align*} 
Observe that the difference with the explicit Euler case lies in the definition of $\uncst$. From~\eqref{neumannfulldiscbackward} one clearly sees that $\unlin$ and $\uncst$ then satisfy again~\eqref{neumannfulldiscforwardinterp} with $\fncst(\bx,t)$ replacing $f_n(\bx)$.

The following estimates holds. 
\begin{lem} 
Consider problem~\eqref{neumannfulldiscbackward} with kernel $\bK$ and data $(\bff,\bg)$ and step-sizes $\tauh > 0$ for all $k$. Assume that $\injn \bK$ satisfies~\ref{assum:Kpos}-\ref{assum:Ksym} and~\ref{assum:KL1}, that $\injn\bg \in L^{\max(p,q)}(\O)$, for some $q \in [1,+\infty]$, and $\sup_{n \in \N}\norm{\injn\bg}_{L^q(\O)} < +\infty$, and that $\fncst \in L^1([0,T];L^{\max(p,q)}(\O))$ and $\sup_{n \in \N}\norm{\fncst}_{L^1([0,T];L^q(\O))} < +\infty$. Then
\[
\text{$\uncst(\cdot,t) \in L^{\max(p,q)}(\O)$}, \forall t \in [0,T], \qandq \sup_{t \in [0,T], n \in \N} \norm{\uncst(\cdot,t)}_{L^q(\O)} < +\infty .
\]
\label{lem:estimatesuncstbackward}
\end{lem}

\bpf{}
Recall from Proposition~\ref{prop:plap}\ref{item4compaccretive} that $J_{\lambda\plap{\injn\bK}{p}}$, $\lambda > 0$, is single-valued on $L^p(\O)$ and nonexpansive on $L^q (\O)$ for all $q \in [1,+\infty]$. Therefore, by induction, we have that for any $k \in [N]$, 
\begin{align*}
\norm{u_n^{k}}_{L^p(\O)} \leq \norm{\injn\bg}_{L^p(\O)} + \sum_{i=0}^k \tau_i \norm{f_n^i}_{L^p(\O)} &\leq \norm{\injn\bg}_{L^p(\O)} + \sum_{i=0}^N \tau_i \norm{f_n^i}_{L^p(\O)} \\
&= \norm{\injn\bg}_{L^p(\O)} + \norm{\fncst}_{L^1([0,T];L^{p}(\O))} .
\end{align*}
Thus $u_n^{k} \in L^p(\O)$, for all $k \in [N]$. In turn, $J_{\tauh\plap{\injn\bK}{p}}(u_n^{k})$ is single-valued for all $k$, and arguing as above, its nonexpansiveness yields
\begin{align*}
\norm{u_n^{k}}_{L^q(\O)} \leq \norm{\injn\bg}_{L^q(\O)} + \norm{\fncst}_{L^1([0,T];L^q(\O))} .
\end{align*}
Taking the supremum over $k$ and $n$ and using the definition of $\uncst$ and the assumptions on $\bg$ and $\bff$, we conclude.
\epf{}

\begin{lem}
Suppose that the assumptions of Lemma~\ref{lem:estimatesuncstbackward} are satisfied with $q=2$ when $p \in ]1,2]$, $q=2(p-1)$ when $p \geq 2$. Assume in addition that $\sup_{n \in \N}\norm{\injn\bK}_{L^{\infty,1}(\O^2)} < +\infty$ and $\sup_{n \in \N}\norm{\fncst}_{\BV([0,T];L^2(\O))} <+\infty$. Then
\begin{equation*}
\sup_{t \in [0,T], n \in \N}\norm{\unlin(\cdot,t)- \uncst(\cdot,t)}_{L^2(\O)} \leq C \tau,
\end{equation*}
where $C$ is a positive constant that does not depend on $(n,N,T)$.
\label{lem:unlincstbackward}
\end{lem}
\bpf{}
For $t \in ]t_{k-1},t_k]$, $k \in \N$, we have
\begin{align}
\norm{\unlin(\cdot,t) - \uncst(\cdot,t)}_{L^2(\O)} 
&= (t_k - t) \norm{\frac{u_n^{k-1} - u_n^{k}}{\tauhm}}_{L^2(\O)} \nonumber\\
&= (t_k - t) \norm{\plap{\injn \bK}{p}u_n^{k} - f_n^k}_{L^2(\O)} \nonumber\\
&\leq \tau \norm{\plap{\injn \bK}{p}u_n^{k} - f_n^k}_{L^2(\O)} \nonumber\\
&= \tau \norm{\plap{\injn \bK}{p}\uncst(\cdot,t_k) - \fncst(\cdot,t_k)}_{L^2(\O)} \nonumber\\
&\leq \tau \pa{\norm{\plap{\injn \bK}{p}\uncst(\cdot,t_k)}_{L^2(\O)} + \sum_{i=1}^k \norm{\fncst(\cdot,t_i) - \fncst(\cdot,t_{i-1})}_{L^2(\O)} + \norm{\fncst(\cdot,0)}_{L^2(\O)}} \nonumber\\
&\leq \tau \pa{\norm{\plap{\injn \bK}{p}\uncst(\cdot,t_k)}_{L^2(\O)} + \Var_{q}(\fncst) + \norm{\fncst(\cdot,0)}_{L^2(\O)}} \nonumber\\
&= \tau \pa{\norm{\plap{\injn \bK}{p}\uncst(\cdot,t_k)}_{L^2(\O)} + \norm{\fncst}_{\BV([0,T];L^2(\O))}} . \label{eq:unlincstbackward}
\end{align}
For $p \in ]1,2]$, we have from~\eqref{eq:plapholder} that
\begin{align*}
\norm{\plap{\injn \bK}{p}\uncst(\cdot,t_k)}_{L^2(\O)} \leq 2^{p/2} C_2^{1/2}\norm{\K}_{L^{\infty,1}(\O^2)}\norm{\uncst(\cdot,t)}_{L^2(\O)}^{p-1} .
\end{align*}
For $p \geq 2$, it is easy to to show with simple arguments as before that
\begin{align*}
\norm{\plap{\injn \bK}{p}\uncst(\cdot,t_k)}_{L^2(\O)} \leq 2^{p-3/2}\norm{\K}_{L^{\infty,1}(\O^2)} \norm{\uncst(\cdot,t)}_{L^{2(p-1)}(\O)}^{p-1} .
\end{align*}
Inserting the last two estimates in~\eqref{eq:unlincstbackward}, taking the supremum over $t$ and $n$ over both sides, and applying Lemma~\ref{lem:estimatesuncstbackward}, we conclude.
\epf

\begin{rem}
As observed in the case of explicit time-discretization the uniform (over $n$) boundedness assumption made in the last two lemmas hold true if $\bg=\projn g$, $\bK=\projn\K$ and $\bff^k = \tauh^{-1}\int_{t_{k-1}}^{t_k} \projn f(\cdot,t)dt$, where $g$, $f$ and $\K$ verify simple assumptions. Indeed, in this case, we have thanks to Lemma~\ref{lem:projinj} that for any $q \in [1,+\infty]$,
\begin{align*}
\sup_{n \in \N}\norm{\injn\bg}_{L^q(\O)} &\leq \norm{g}_{L^q(\O)}, \quad
\sup_{n \in \N}\norm{\injn\bK}_{L^{\infty,q}(\O^2)} \leq \norm{\K}_{L^{\infty,q}(\O^2)}, \\ 
\sup_{n \in \N}\norm{\fncst}_{L^1([0,T];L^q(\O))} &\leq \norm{f}_{L^1([0,T];L^q(\O))} \qandq \sup_{n \in \N}\norm{\fncst}_{\BV([0,T];L^q(\O))} \leq \norm{f}_{\BV([0,T];L^q(\O))} .
\end{align*}
In fact, the condition $f \in \BV([0,T];L^q(\O))$ is sufficient to ensure that 
\[
\sup_{n \in \N}\norm{\fncst}_{L^1([0,T];L^q(\O))} < +\infty \qandq \sup_{n \in \N}\norm{\fncst}_{\BV([0,T];L^q(\O))} < +\infty .
\] 
Indeed, arguing as in~\cite[Lemma~A.1]{BrezisBook}, this conditions implies $f \in L^\infty([0,T];L^q(\O))$. In turn, using Lemma~\ref{lem:projinj}, we get
\begin{align*}
\norm{\fncst}_{L^1([0,T];L^q(\O))} \leq \norm{f}_{L^1([0,T];L^q(\O))} 
&\leq \norm{f}_{L^\infty([0,T];L^q(\O))} \\
&\leq \norm{f(\cdot,0)}_{L^q(\O)} + \Var_{q}(f) = \norm{f}_{\BV([0,T];L^q(\O))} .
\end{align*}
\label{rem:estimatesuncstbackward}
\end{rem}

We are now in position to state the error bound for the fully discrete problem with backward/implicit Euler time discretization.
\begin{theo}
Suppose that $p \in ]1, +\infty[$. Let $u$ be a solution of~\eqref{neumann} with kernel $\K$ and data $(f,g)$, and $\acc{\bu^k}_{k \in [N]}$ is the sequence generated by~\eqref{neumannfulldiscbackward} with $\bK=\projn \K$, $\bg=\projn g$, $\bff^k = \tauh^{-1}\int_{t_{k-1}}^{t_k} \projn f(\cdot,t)dt$. Assume that $\K$ satisfies~\ref{assum:Kpos},~\ref{assum:Ksym} and $\K \in L^{\infty,2}(\O^2)$, and that $f,g$ satisfy either one of the conditions~\ref{assum:maina},~\ref{assum:mainb} or~\ref{assum:mainc} in Theorem~\ref{thm:maincontinuous}, and that $f \in \BV([0,T];L^2(\O))$. Then, the following hold.
\begin{enumerate}[label=(\roman*)]
\item $u$ is the unique solution of~\eqref{neumann}, $\acc{\bu^k}_{k \in [N]}$ is uniquely defined and $\acc{\norm{\injn\bu^k}_{L^2(\O)}}_{k \in [N]}$ is bounded uniformly in $n$.
\item We have the error estimate
\begin{multline}
\sup_{k \in [N], t \in ]t_{k-1},t_k]}\norm{\injn\bu^{k} - u(\cdot,t)}_{L^2(\O)} 
\leq \exp\pa{T/2}\Bigg(\norm{\injn\projn g - g}_{L^2(\O)} + \norm{\fncst - f}_{L^1([0,T];L^2(\O))} \\
+ CT^{1/2}
\begin{cases}
\tau^{1/(3-p)} + \norm{\injn\projn \K - \K}_{L^{\infty,2}(\O^2)}	& \text{under~\ref{assum:maina}} \\
\tau^{p/(2p-1)} + \norm{\injn\projn \K - \K}_{L^{\infty,2}(\O^2)}	& \text{under~\ref{assum:mainb}}\\
\tau^{1/(3-p)} + \norm{\injn\projn \K - \K}_{L^2(\O^2)}		& \text{under~\ref{assum:mainc} when $p \in ]1,2]$} \\
\tau + \norm{\injn\projn \K - \K}_{L^2(\O^2)}				& \text{under~\ref{assum:mainc} when $ p \geq 2 $} .
\end{cases}
\Bigg) ,
\label{eq:fulldiscboundbackwardcst}
\end{multline}
for $\tau$ sufficiently small, where $C$ is positive constant that depends only on $p$, $g$, $f$ and $\K$.
\item If, moreover, $g \in L^{\infty}(\O) \cap \Lip(s,L^2(\O))$, $\K \in \Lip(s,L^2(\O^2))$, and $f \in L^1([0,T];L^{\infty}(\O)) \cap \Lip(s,L^2(\O \times [0,T]))$ then
\begin{multline}
\hspace*{-1cm}
\sup_{k \in [N], t \in ]t_{k-1},t_k]}\norm{\injn\bu^{k} - u(\cdot,t)}_{L^2(\O)}\leq C\exp(T/2)\Bigg((1+T^{1/2})\deltan^{s} \\
+ T^{1/2}\Bigg(
\begin{cases}
\tau^{\min(s,1/(3-p))}		& \text {when $p \in ]1,2]$} \\
\tau^s				& \text{ when $ p \geq 2 $} 
\end{cases}\Bigg)\Bigg).
\label{eq:fulldiscretebackwardratecst}
\end{multline}
for $\tau$ sufficiently small, where $C$ is positive constant that depends only on $p$, $g$, $f$, $\K$ and $s$. The term $\tau^s$ in the dependence on $\tau$ disappears when $f$ is time-independent.
\end{enumerate}
\label{thm:mainfulldiscretebackward}
\end{theo}

\bpf{}
In the proof, $C$ is any positive constant that may depend solely on $p$, $g$, $f$, $\K$ and/or $s$, and that may be different at each line.
\begin{enumerate}[label=(\roman*)]
\item Existence and uniqueness of $u$ were proved in Theorem~\ref{thm:maincontinuous}\ref{item:maincontinuouswellposed}. Well-posedness of the sequence $\acc{\bu^k}_{k \in [N]}$ is a consequence of Lemma~\ref{lem:estimatesuncstbackward} and Remark~\ref{rem:estimatesuncstbackward}.

\item For $p \in ]1,2]$, the proof of the error bound is exactly the same as that of~\eqref{eq:fulldiscboundforwardcst} in Theorem~\ref{thm:mainfulldiscreteforward} using the modified definition of $\uncst$ and that now $f$ is time-dependent, and thus we replace $f_n$ there by $\fncst$. We also denote $g_n = \injn\projn g$ and $\Kn = \injn\projn\K$.

For the $p \geq 2$, the argument is also similar, and the main change consists in bounding appropriately the third term in~\eqref{eq:derivfullforward}. We then invoke Lemma~\ref{lem:unlincstbackward} to show that
\begin{align*}
\abs{\int_{\O} \pa{\plap{\Kn}{p} \uncst(\bx,t) - \plap{\Kn}{p} u(\bx,t)}\pa{\unlin(\bx,t) - \uncst(\bx,t)} d\bx} \leq C\norm{\plap{\Kn}{p} \uncst(\cdot,t)-\plap{\Kn}{p}u(\cdot,t)}_{L^2(\O)}\tau ,
\end{align*}
where $C$ is indeed a finite constant owing to the assumption on $f$ and Remark~\ref{rem:estimatesuncstbackward}. We now use Lemma~\ref{lem:psiprop}\ref{item:psicont} to get the bound
\begin{align}
&\norm{\plap{\Kn}{p} \uncst(\cdot,t)-\plap{\Kn}{p}u(\cdot,t)}_{L^2(\O)}^2 \nonumber\\ 
&= \int_{\O} \abs{\int_{\O}\Kn(\bx,\by) \pa{\Psi(\uncst(\by,t)-\uncst(\bx,t))-\Psi(u(\by,t)-u(\bx,t))} d\by}^2 d\bx \nonumber\\
&\leq \int_{\O} \pa{\int_{\O}\Kn(\bx,\by) \abs{\xicst(\by,t)-\xicst(\bx,t))}\pa{|\uncst(\by,t)-\uncst(\bx,t)|+|u(\by,t)-u(\bx,t)|}^{p-2} d\by}^2 d\bx .
\label{eq:plapcontpgt2}
\end{align}
For case~\ref{assum:mainc}, we infer from Lemma~\ref{lem:estimatesuncstbackward} (with $q=+\infty$) and Lemma~\ref{lem:projinj} that
\begin{align}
&\norm{\plap{\Kn}{p} \uncst(\cdot,t)-\plap{\Kn}{p}u(\cdot,t)}_{L^2(\O)}^2 \nonumber\\
&\leq \pa{4\pa{\norm{g}_{L^{\infty}(\O)}+\norm{f}_{L^1([0,T];L^{\infty}(\O))}}}^{2(p-2)} 
\int_{\O} \pa{\int_{\O}\Kn(\bx,\by) \abs{\xicst(\by,t)-\xicst(\bx,t))} d\by}^2 d\bx \nonumber\\
&\leq \pa{4\pa{\norm{g}_{L^{\infty}(\O)}\norm{f}_{L^1([0,T];L^{\infty}(\O))}}}^{2(p-2)} \norm{\K}_{L^{\infty,2}(\O^2)} \int_{\O^2}\Kn(\bx,\by) \abs{\xicst(\by,t)-\xicst(\bx,t))} ^2 d\bx d\by \nonumber\\
&= 4\pa{4\pa{\norm{g}_{L^{\infty}(\O)}\norm{f}_{L^1([0,T];L^{\infty}(\O))}}}^{2(p-2)} \norm{\K}_{L^{\infty,2}(\O^2)} \int_{\O^2}\Kn(\bx,\by) \abs{\xicst(\bx,t)} ^2 d\bx d\by \nonumber\\
&\leq 4\pa{4\pa{\norm{g}_{L^{\infty}(\O)}\norm{f}_{L^1([0,T];L^{\infty}(\O))}}}^{2(p-2)} \norm{\K}_{L^{\infty,2}(\O^2)}^2 \norm{\xicst(\cdot,t)}_{L^2(\O)}^2 .
\label{eq:plapcontpgt2bnd}
\end{align}
It then follows by Cauchy-Schwartz inequality that
\begin{align*}
& \abs{\int_{\O} \pa{\plap{\Kn}{p} \uncst(\bx,t) - \plap{\Kn}{p} u(\bx,t)}\pa{\unlin(\bx,t) - \uncst(\bx,t)} d\bx} \\
&\leq C\norm{\xicst(\cdot,t)}_{L^2(\O)}\tau \\
&\leq C\pa{\norm{\xilin(\cdot,t)}_{L^2(\O)}\tau + \tau^2} \\
&\leq \frac{1}{6}\norm{\xilin(\cdot,t)}_{L^2(\O)}^2 + C\tau^2 .
\end{align*}
Inserting this in~\eqref{eq:derivfullforward}, using again Young inequality for the last term, we have shown that when $p \geq 2$ and~\ref{assum:mainc} holds,
\begin{equation*}
\frac{\partial}{\partial t} \norm{\xilin(\cdot,t)}_{L^2(\O)}^2 
\leq \norm{\xilin(\cdot,t)}_{L^2(\O)}^2 + C\pa{\tau^2 + \norm{\fncst(\cdot,t) - f(\cdot,t)}_{L^2(\O)}^2 + \norm{\Kn - \K}_{L^2(\O^2)}^2} .
\end{equation*}
Using the Gronwall's lemma, taking the square-root and using~\eqref{eq:bnduncstulin}, we get the error bound in this case.

It remains to consider the case~\ref{assum:mainb}, when $p\geq2$. For this, we embark from~\eqref{eq:plapcontpgt2}, and use the continuity of $\Psi$ in Lemma~\ref{lem:psiprop}~\ref{item:psimonotone} (see \eqref{eq:psigencont}) with $\alpha=1/p$. Combining this with Jensen and H\"older inequalities, we get
\begin{align*}
&\norm{\plap{\Kn}{p} \uncst(\cdot,t)-\plap{\Kn}{p}u(\cdot,t)}_{L^2(\O)}^2 \\
&\leq \norm{\K}_{L^{\infty,1}(\O^2)} \int_{\O^2}\pa{\Kn(\bx,\by) \abs{\xicst(\by,t)-\xicst(\bx,t))}^{2/p}}\\
&\qquad \pa{|\uncst(\by,t)-\uncst(\bx,t)|+|u(\by,t)-u(\bx,t)|}^{2(p-1)-2/p} d\bx d\by \\
&\leq \norm{\K}_{L^{\infty,1}(\O^2)} \int_{\O^2}\pa{\Kn(\bx,\by) \abs{\xicst(\by,t)-\xicst(\bx,t))}^{2}}^{1/p}\\
&\qquad  \pa{\pa{\Kn(\bx,\by)}^{(p-1)/p}\pa{|\uncst(\by,t)-\uncst(\bx,t)|+|u(\by,t)-u(\bx,t)|}^{2(p-1)-2/p}} d\bx d\by \\
&\leq \norm{\K}_{L^{\infty,1}(\O^2)} \pa{\int_{\O^2}\Kn(\bx,\by) \abs{\xicst(\by,t)-\xicst(\bx,t))}^2 d\bx d\by}^{1/p} \\ 
&\qquad \pa{\int_{\O^2}\Kn(\bx,\by)\pa{|\uncst(\by,t)-\uncst(\bx,t)|+|u(\by,t)-u(\bx,t)|}^{2p-2/(p-1)} d\bx d\by}^{(p-1)/p} \\
&\leq \norm{\K}_{L^{\infty,1}(\O^2)} \pa{4\int_{\O^2}\Kn(\bx,\by) \abs{\xicst(\bx,t))}^2 d\bx d\by}^{1/p} \\ 
&\qquad \pa{2^{2p-2/(p-1)}\int_{\O^2}\Kn(\bx,\by)\pa{\abs{\uncst(\bx,t)}+\abs{u(\bx,t)}}^{2p-2/(p-1)} d\bx d\by}^{(p-1)/p} .\\
&\leq 4\norm{\K}^{2}_{L^{\infty,1}(\O^2)} \pa{\int_{\O^2}\abs{\xicst(\bx,t))}^2 d\bx d\by}^{1/p} \\ 
&\qquad \pa{\int_{\O^2}\pa{\abs{\uncst(\bx,t)}+\abs{u(\bx,t)}}^{2p-2/(p-1)} d\bx d\by}^{(p-1)/p} .
\end{align*}
Observe that $L^{2p-2/(p-1)}(\O)\subset L^{2(p-1)}(\O)$, hence by H\"older inequality and Lemma~\ref{lem:estimatesuncstbackward} with $q=2(p-1)$ and Lemma~\ref{lem:projinj}, the last term in the above display can be bounded as 
\begin{align*}
&\pa{\int_{\O^2}\pa{\abs{\uncst(\bx,t)}+\abs{u(\bx,t)}}^{2p-2/(p-1)} d\bx d\by}^{(p-1)/p}\\
&\leq \norm{\abs{\uncst(\bx,t)}+\abs{u(\bx,t)}}_{L^{2(p-1)}(\O)}^{2(p-1)-2/p} \\
&\leq \pa{\norm{g}_{L^{2(p-1)}(\O)}+\norm{f}_{L^1([0,T];L^{2(p-1)}(\O))}}^{2(p-1)-2/p} .
\end{align*}
We then arrive at 
\begin{align*}
\norm{\plap{\Kn}{p} \uncst(\cdot,t)-\plap{\Kn}{p}u(\cdot,t)}_{L^2(\O)}^2 \leq C \norm{\K}^2_{L^{\infty,1}(\O^2)} \norm{\xicst}_{L^2(\O)}^{2/p}.
\end{align*}
Hence 
\begin{align*}
&\abs{\int_{\O} \pa{\plap{\Kn}{p} \uncst(\bx,t) - \plap{\Kn}{p} u(\bx,t)}\pa{\unlin(\bx,t) - \uncst(\bx,t)} d\bx} \\
&\leq C\norm{\xicst(\cdot,t)}_{L^2(\O)}^{1/p}\tau \\
&\leq C\pa{\norm{\xilin(\cdot,t)}_{L^2(\O)}^{1/p}\tau+ \tau^{(p+1)/p}}\\
&\leq \frac{1}{6}\norm{\xilin(\cdot,t)}_{L^2(\O)}^2 + C(\tau^{2p/(2p-1)}+ \tau^{(p+1)/p}) .
\end{align*}
Inserting this into~\eqref{eq:derivfullforward}, using again Young inequality for the last term,
\begin{multline*}
\frac{\partial}{\partial t} \norm{\xilin(\cdot,t)}_{L^2(\O)}^2 
\leq \norm{\xilin(\cdot,t)}_{L^2(\O)}^2 + C\Bigg(\tau^{2p/(2p-1)}+ \tau^{(p+1)/p} + \norm{\fncst(\cdot,t) - f(\cdot,t)}_{L^2(\O)}^2\\
+ \norm{\Kn - \K}_{L^2(\O^2)}^2 \Bigg).
\end{multline*}
Hence, using the Gronwall's lemma, taking the square-root and using~\eqref{eq:bnduncstulin}, we get the error bound in this case, after observing that the dependence on $\tau$ scales as $O(\tau^{p/(2p-1)})$ for $\tau$ sufficiently small (or $N$ large enough) since $1/2 < p/(2p-1) \leq (p+1)/(2p)$ for $p \geq 2$.

\item Plug~\eqref{eq:lipspaceapprox} into~\eqref{eq:fulldiscboundbackwardcst} after observing that
\[
\hspace*{-1cm}
\norm{\fncst - f}_{L^1([0,T];L^2(\O))} \leq T^{1/2} \norm{\fncst - f}_{L^2([0,T];L^2(\O))} = T^{1/2}\norm{\fncst - f}_{L^2(\O \times [0,T])} \leq C T^{1/2}\max(\tau^s,\deltan^s) .
\] 
For the scaling in $\tau$, we use that $s \in ]0,1]$. 
\end{enumerate}
\epf

Another way to derive error bounds for~\eqref{neumannfulldiscbackward} is as follows. To lighten notation, denote $g_n = \injn\projn g$, $f_n(\cdot,t) = \injn\projn f(\cdot,t)$ for $t \in [0,T]$, and $\Kn = \injn\projn \bK$. Let $u_n$ be a solution to~\eqref{neumann} with data $(f_n,g_n)$ and kernel $\Kn$. Under the assumptions of Theorem~\ref{thm:mainfulldiscretebackward} on $(f,g,\K)$, $u_n$ is unique. Then one has
\[
\norm{\unlin - u}_{C([0,T];L^2(\O))} \leq \norm{\unlin - u_n}_{C([0,T];L^2(\O))} + \norm{u_n - u}_{C([0,T];L^2(\O))} .
\]
Theorem~\ref{thm:maincontinuous} provides a bound on the last term of the right-hand side in the above display, which captures the space-discretization error.
Bounds for the first term, which corresponds to the time-discretization error, were derived in $C([0,T];L^p(\O))$ by Crandall and Liggett in their seminal paper~\cite{CrandallLiggett71} for constant time step-size and $f = 0$, and then extended to non-uniform time partitions in~\cite{kobayashi}, see also~\cite{Savare06}. More precisely, using~\cite[Theorem~1]{Savare06} and the fact that $\unlin(\cdot,0) = u_n(\cdot,0) = g_n$, the following bound holds
\begin{multline*}
\norm{\unlin - u_n}_{C([0,T];L^p(\O))} \leq \norm{\fncst - f_n}_{L^1([0,T];L^p(\O))} + 2T^{1/2}\pa{\norm{f_n^1 - \plap{\Kn}{p}g_n}_{L^p(\O)}+\Var_{p}(\fncst)}\tau^{1/2}.
\end{multline*}
The first term can be bounded as follows (for constant step-size to simplify)
\begin{align*}
\norm{\fncst - f_n}_{L^1([0,T];L^p(\O))} 
&= \sum_{k=1}^N \int_{t_{k-1}}^{t_k} \norm{\tauh^{-1}\int_{t_{k-1}}^{t_k} f_n(\cdot,s)ds - f_n(\cdot,t)}_{L^p(\O)} dt \\
&\leq \tau^{-1}\sum_{k=1}^N \int_{t_{k-1}}^{t_k}\int_{t_{k-1}}^{t_k} \norm{f(\cdot,s) - f(\cdot,t)}_{L^p(\O)} dsdt \\
&\leq \tau^{-1} \int_{-\tau}^{\tau} \pa{\int_{0}^T \norm{f(\cdot,t+s) - f(\cdot,t)}_{L^p(\O)} dt}ds \\
&\leq \tau^{-1} \int_{-\tau}^{\tau} s \Var_{p}(f) ds = \tau \Var_{p}(f) ,
\end{align*} 
where we used Lemma~\ref{lem:projinj} in the first inequality and~\cite[Lemma~A.1]{BrezisBook} in the last one. Overall, this shows that the time discretization error $\norm{\unlin - u_n}_{C([0,T];L^p(\O))}$ scales as $O\pa{\pa{T\tau}^{1/2}}$ for $\tau$ sufficiently small. The rate $O(\tau^{1/2})$ is known to be optimal for general accretive operators in Banach spaces (see~\cite{Savare06}). In turn, by standard comparisons of $L^q(\O)$ norms (assuming that~\ref{assum:mainc} holds so that boundedness of $\unlin$ and $u_n$ is in force), this strategy gives us a bound which scales as 
\[
\norm{\unlin - u_n}_{C([0,T];L^p(\O))} = 
\begin{cases}
O\pa{\tau^{1/2}} & p \geq 2, \\
O\pa{\tau^{p/4}} & p \in ]1,2] .
\end{cases}
\]
This is strictly worse than the rates in $\tau$ obtained from~\eqref{eq:fulldiscboundbackwardcst}. There is however no contradiction in this and the reason is that the strategy outlined above is too general and does not exploit all properties of the operator $\plap{\K}{p}$ among which its continuity that was a key to derive better rates in $\tau$. In this sense, our present results are optimal. We also remark that our rates are consistent with those in~\cite{Hafienne_FE_Nonlocal_plaplacian_evolu} for $p \geq 2$.


\section{Application to random graph sequences}
\label{sec:applicationgraphs}
In this section, we study continuum limits of fully discrete problems on the random graph model of Definition~\ref{def:randomgraph} with backward/implicit Euler time discretization. Explicit discretization can also be treated following our results in Section~\ref{subsec:forwardeuler}, but we will not elaborate further on it for the sake of brevity.

Recall the notations in Section~\ref{graphlimits}, in which case we now set $\O = [0,1]$. Recall also the the construction of the random graph model in Definition~\ref{def:randomgraph} where each edge $(i,j)$ is independently set to $1$ with probability~\eqref{eq:sparsegraphmodelavg}. This entails that the random matrix $\bLam$ is symmetric. However, it is worth emphasizing that the entries of $\bLam$ are not independent, but only the entries in each row are mutually independent\footnote{This feature was already used in the proof of Lemma~\ref{lem:Linf1norm}}. This observation will be instrumental in deducing our error bound.

We consider the fully discrete on $\K$-random graphs $\bG(n,\K,\rho_n)$ with backward Euler time scheme 
\begin{equation}\tag{\textrm{${\mathcal{P}}^{\text{TDB},\bG}_{p}$}}
\begin{cases}
\displaystyle{\frac{\bu^{k} - \bu^{k-1}}{\tauhm}} = \frac{1}{\rho_n n}\sum\limits_{j: (i,j) \in E(\bG(n,\K,\rho_n))} \Psi(\bu_{\bj} - \bu_{\bi}) + \bff^k, & k \in [N], \\
\bu^0 = \bg ,
\end{cases}
\label{neumannfulldiscbackwardgraphs}
\end{equation}
where $\bu^k, \bff^k \in \R^{n}$. It is important to keep in mind that, since $\bG(n,\K,\rho_n)$ is a random variable taking values in the set of simple graphs, the evolution problem~\eqref{neumannfulldiscbackwardgraphs} must be understood in this sense. Observe that the normalization in~\eqref{neumannfulldiscbackwardgraphs} by $\rho_n n$ corresponds to the average degree (see Section~\ref{subsec:randomgraph} for details).

Problem \eqref{neumannfulldiscbackwardgraphs} can be equivalently written as
\begin{equation*}
\begin{cases}
\displaystyle{\frac{\bu^{k} - \bu^{k-1}}{\tauhm}} = -\plapdisc{\bLam}{p}\bu^{k} + \bff^k, & k \in [N], \\
\bu^0 = \bg .
\end{cases}
\end{equation*}

We define the time-space continuum extensions $\unlin$ and $\uncst$ and as in Section~\ref{subsec:backwardeuler}. One then sees that they satisfy
\begin{equation}
\begin{cases}
\frac{\partial}{\partial t} \unlin(x,t)=-\plap{\injn\bLam}{p}\uncst(x,t) + \fncst(x,t), & (x,t) \in \O \times ]0,T], \\
\unlin(x,0) = \injn\bg(x), \quad x \in \O .
\end{cases}
\label{neumannfulldiscbackwardinterpgraphs}
\end{equation}
Toward our goal of establishing error bounds, we define $\bv$ as the solution of the fully discrete problem~\eqref{neumannfulldiscbackward} with data $(\bff,\bg)$ and discrete kernel $\wedg$. Its time-space continuum extensions, $\vnlin$ and $\vncst$, defined similarly as above, fulfill
\begin{equation}
\begin{cases}
\frac{\partial}{\partial t} \vnlin(x,t)=-\plap{\injn\wedg}{p}\vncst(x,t) + \fncst(x,t), & (x,t) \in \O \times ]0,T], \\
\vnlin(x,0) = \injn\bg(x), \quad x \in \O .
\end{cases}
\label{neumannintermediatebackwardinterpgraphs}
\end{equation}
We have 
\begin{equation}\label{eq:errorbnddecomp}
\norm{\unlin - u}_{C([0,T];L^2(\O))} \leq \norm{\unlin - \vnlin}_{C([0,T];L^2(\O))} + \norm{\vnlin - u}_{C([0,T];L^2(\O))} .
\end{equation}
This bound is composed of two terms: the first one captures the error of random sampling, and the second that of (space and time) discretization. We start by bounding the first term by comparing~\eqref{neumannfulldiscbackwardinterpgraphs} and~\eqref{neumannintermediatebackwardinterpgraphs}.

\begin{lem}
Assume that $(\bff^k,\bg,\bK,f,g,\K)$ verify the assumptions of Theorem~\ref{thm:mainfulldiscretebackward}. Assume also that $\rho_n \to 0$ and $n\rho_n = \omega\pa{\pa{\log n}^\gamma}$ for some $\gamma > 1$. Then, for any $\beta \in ]0,1[$,
\begin{equation}\label{eq:intermediatebndfulldiscretebackwardgraph}
\norm{\unlin - \vnlin}_{C([0,T];L^2(\O))}
\leq C\exp\pa{T/2}T^{1/2} \pa{(\rho_n n)^{-\beta/2} +
\begin{cases}
\tau^{1/(3-p)} 	& p \in ]1,2], \\
\tau 		& p \geq 2 .
\end{cases}} .
\end{equation}
with probability at least $1-(\rho_n n)^{-(1-\beta)}$. In particular,
\begin{equation}\label{eq:intermediatebndfulldiscretebackwardgraph2}
\norm{\unlin - \vnlin}_{C([0,T];L^2(\O))}
\leq C\exp\pa{T/2}T^{1/2} \pa{o\pa{\pa{\log n}^{-\gamma\beta/2}} +
\begin{cases}
\tau^{1/(3-p)} 	& p \in ]1,2], \\
\tau 		& p \geq 2 .
\end{cases}} .
\end{equation}
with probability at least $1-o\pa{\pa{\log n}^{-\gamma(1-\beta)}}$.
\label{lem:intermediatebndfulldiscretebackwardgraphs}
\end{lem}

\bpf{}
Denote $\xilin(x,t) = \vnlin(x,t) - \unlin(x,t)$, $\xicst(x,t) = \vncst(x,t) - \uncst(x,t)$, $g_n = \injn\projn g$, $\Wedg = \injn\wedg$ and $\Lambdan = \injn\bLam$. We thus have from~\eqref{neumannfulldiscbackwardinterpgraphs} and~\eqref{neumannintermediatebackwardinterpgraphs} that a.e.
\begin{equation*}
\begin{split}
\frac{\partial \xilin(x,t)}{\partial t } 
&= - \pa{\plap{\Wedg}{p}(\vncst(x,t)) - \plap{\Lambdan}{p}(\uncst(x,t))} \\
&= - \pa{\plap{\Wedg}{p}(\vncst(x,t)) - \plap{\Lambdan}{p}(\vncst(x,t))} - \pa{\plap{\Lambdan}{p}(\vncst(x,t)) - \plap{\Lambdan}{p}(\uncst(x,t))} .
\end{split}
\end{equation*}
Multiplying both sides by $\xilin(x,t)$, integrating and rearranging the terms, we get
\begin{equation}
\begin{split}
&\frac{1}{2} \frac{\partial}{\partial t} \norm{\xilin(\cdot,t)}_{L^2(\O)}^2 
= - \int_{\O} \pa{\plap{\Lambdan}{p} \vncst(x,t) - \plap{\Lambdan}{p} \uncst(x,t)}(\vncst(x,t) - \uncst(x,t)) dx \\
& - \int_{\O} \pa{\plap{\Wedg}{p} \vncst(x,t) - \plap{\Lambdan}{p} \vncst(x,t)}\xilin(x,t) dx \\
& - \int_{\O} \pa{\plap{\Lambdan}{p} \vncst(x,t) - \plap{\Lambdan}{p} \uncst(x,t)}\pa{(\vnlin(x,t) - \vncst(x,t)) - (\unlin(x,t) - \uncst(x,t))} dx .
 \end{split}
\label{eq:derivfullbackwardgraph}
\end{equation}
Under our condition on $n\rho_n$, Lemma~\ref{lem:Linf1norm} tells us that with probability 1,
\[
\norm{\Lambdan}_{L^{\infty,1}(\O^2)} = \norm{\Wedg}_{L^{\infty,1}(\O^2)} + o(1) \leq \norm{\injn\projn\K}_{L^{\infty,1}(\O^2)} + o(1) \leq \norm{\K}_{L^{\infty,1}(\O^2)} + o(1) ,
\] 
so in particular $\norm{\Lambdan}_{L^{\infty,1}(\O^2)}$ is uniformly bounded with probability 1. $\Lambdan$ is also positive and symmetric. Since $g \in L^{q}(\O)$ and $f\in L^{1}([0,T];L^{q}(\O))\cap\BV([0,T];L^{2}(\O))$, $q \in \{2,2(p-1),+\infty\}$, the conclusions of Lemma~\ref{lem:estimatesuncstbackward} and Lemma~\ref{lem:unlincstbackward} remain true which shows that with probability 1,
\[
\sup_{t \in [0,T], n \in \N}\norm{\uncst(\cdot,t)}_{L^{q}(\O)} < +\infty \qandq \sup_{t \in [0,T], n \in \N}\norm{\unlin(\cdot,t)- \uncst(\cdot,t)}_{L^2(\O)} \leq C \tau .
\] 
The same claim holds for $\vnlin$ and $\vncst$ since $\norm{\Wedg}_{L^{\infty,1}(\O^2)} \leq \norm{\K}_{L^{\infty,1}(\O^2)} < +\infty$ and $\Wedg$ is positive and symmetric, i.e. $\Wedg$ obeys~\ref{assum:Kpos}-\ref{assum:KL1}. Thus Proposition~\ref{prop:plap}\ref{itemmonotonep} entails that the first term on the right-hand side of~\eqref{eq:derivfullbackwardgraph} is nonpositive with probability 1. Let us now bound the second term. Denote the random variables $\bZ_i \eqdef \frac{1}{n}\sum_{j \in [n]}\bpa{\bLam_{ij} - \wedg_{ij}}\Psi(\bv_{j} - \bv_{i})$. By Cauchy-Schwartz inequality, we have
\begin{align*}
\abs{\int_{\O} \pa{\plap{\Lambdan}{p} \uncst(x,t) - \plap{\Wedg}{p} \uncst(x,t)}\xilin(x,t) dx} 
\leq C \norm{\injn\bZ}_{L^{2}(\O)} \norm{\xilin(\cdot,t)}_{L^{2}(\O)} .
\end{align*}
For the last term in~\eqref{eq:derivfullbackwardgraph}, we argue as in the proof of Theorem~\ref{thm:mainfulldiscretebackward} to show that, with probability 1, 
\begin{multline*}
\abs{\int_{\O} \pa{\plap{\Lambdan}{p} \uncst(x,t) - \plap{\Lambdan}{p} \vncst(x,t)}\pa{(\unlin(x,t) - \uncst(x,t)) - (\vnlin(x,t) - \vncst(x,t))} dx} \\
\leq C
\begin{cases}
{\norm{\xilin(\cdot,t)}_{L^2(\O)}^{p-1}\tau + \tau^p} & p \in ]1,2], \\
{\norm{\xilin(\cdot,t)}_{L^2(\O)}\tau + \tau^2} & p \geq 2 .
\end{cases}
\end{multline*}
Collecting all these bounds, after using Young inequality, we have shown that (again with probability 1),
\begin{equation*}
\frac{\partial}{\partial t} \norm{\xilin(\cdot,t)}_{L^2(\O)}^2 
\leq 
\norm{\xilin(\cdot,t)}_{L^2(\O)}^2 +
C \pa{\norm{\injn\bZ}_{L^{2}(\O)}^2 + 
\begin{cases}
\tau^{2/(3-p)} + \tau^p & p \in ]1,2], \\
\tau^2 & p \geq 2 .
\end{cases}}
\end{equation*}
Using the Gronwall's lemma and taking the square-root, we get for $\tau$ sufficiently small
\begin{equation}\label{eq:unvngraphgronwall}
\norm{\unlin - \vnlin}_{C([0,T];L^2(\O))}
\leq C\exp\pa{T/2}T^{1/2} \pa{\norm{\injn\bZ}_{L^{2}(\O)} +
\begin{cases}
\tau^{1/(3-p)} 	& p \in ]1,2], \\
\tau 		& p \geq 2 .
\end{cases}} .
\end{equation}
It remains to bound the random variable $\norm{\injn\bZ}_{L^{2}(\O)}$. For this purpose, we have by Markov inequality that for $\varepsilon > 0$
\begin{align*}
\P\pa{\norm{\injn\bZ}_{L^{2}(\O)} \geq \varepsilon} = \P\pa{n^{-1}\sum_{i}\bZ_i^2 \geq \varepsilon^2} \leq \varepsilon^{-2} n^{-1}\sum_{i}\EE\pa{\bZ_i^2} .
\end{align*}
By independence of $\pa{\bLam_{ij}}_{j \in [n]}$, for each $i \in [n]$, we get
\begin{align*}
\EE\pa{\bZ_i^2} 
&= (\rho_n n)^{-2} \sum_{j \in [n]} \VV\pa{\rho_n\bLam_{ij}}\pa{\Psi(\bv_{j} - \bv_{i})}^2 = (\rho_n n)^{-2} \sum_{j\in [n]} \rho_n\wedg_{ij}(1-\rho_n\wedg_{ij})\pa{\Psi(\bv_{j} - \bv_{i})}^2 \\
&\leq (\rho_n n^2)^{-1} \sum_{j \in [n]} \wedg_{ij}\abs{\bv_{j} - \bv_{i}}^{2(p-1)}.
\end{align*}
In turn,
\begin{align*}
\P\pa{\norm{\injn\bZ}_{L^{2}(\O)} \geq \varepsilon}& \leq (\varepsilon^{2}\rho_n n)^{-1}\frac{1}{n^2} \sum_{i,j \in [n]} \wedg_{ij}\abs{\bv_{j} - \bv_{i}}^{2(p-1)}\\
& = (\varepsilon^{2}\rho_n n)^{-1}\int_{\O^2}\Wedg(\bx,\by)\abs{\vncst (\by) - \vncst (\bx)}^{2(p-1)}d\by d\bx .
\end{align*}
If the condition~\ref{assum:maina} holds, then by the symmetry of the kernel, Jensen inequality and H\"older inequality, one gets
\begin{align*}
\int_{\O^2}\Wedg(\bx,\by)\abs{\vncst (\by) - \vncst (\bx)}^{2(p-1)}d\by d\bx 
&\leq 4\int_{\O^2}\Wedg(\bx,\by)\abs{\vncst (\bx)}^{2(p-1)}d\by d\bx\\
&\leq 4\norm{\Wedg}_{L^{\infty,1}(\O^2)}\int_{\O}\abs{\vncst (\bx)}^{2(p-1)}d\bx\\
&\leq 4\norm{\Wedg}_{L^{\infty,1}(\O^2)}\norm{\vncst}^{2(p-1)}_{L^{2}(\O)}.
\end{align*}
Under the condition~\ref{assum:mainb}, by the symmetry of the kernel and Jensen inequality again, we have
\begin{align*}
\int_{\O^2}\Wedg(\bx,\by)\abs{\vncst (\by) - \vncst (\bx)}^{2(p-1)}d\by d\bx &\leq 2^{2(p-1)}\int_{\O^2}\Wedg(\bx,\by)\abs{\vncst (\bx)}^{2(p-1)}d\by d\bx\\
&\leq 2^{2(p-1)}\norm{\Wedg}_{L^{\infty,1}(\O^2)}\norm{\vncst}^{2(p-1)}_{L^{2(p-1)}(\O)}.
\end{align*}
Similarly, under condition~\ref{assum:mainc}, we have
\begin{align*}
\int_{\O^2}\Wedg(\bx,\by)\abs{\vncst (\by) - \vncst (\bx)}^{2(p-1)}d\by d\bx 
&\leq 2^{2(p-1)}\norm{\vncst}^{2(p-1)}_{L^{\infty}(\O)} \norm{\Wedg}_{L^1(\O^2)} \\
&\leq 2^{2(p-1)}\norm{\Wedg}_{L^{\infty,1}(\O^2)} \norm{\vncst}^{2(p-1)}_{L^{\infty}(\O)} .
\end{align*}
Since $\norm{\Wedg}_{L^{\infty,1}(\O^2)}\leq \norm{\K}_{L^{\infty,1}(\O^2)}$ (see~\eqref{proppvnbi} in Lemma~\ref{lem:projinj}), we have 
\begin{align*}
\P\pa{\norm{\injn\bZ}_{L^{2}(\O)} \geq \varepsilon}& \leq C (\varepsilon^{2}\rho_n n)^{-1} \norm{\K}_{L^{\infty,1}(\O^2)} ,
\end{align*}
where 
\begin{align*}
C=
\begin{cases}
4\sup_n\norm{\vncst}^{2(p-1)}_{L^{2}(\O)}, 				& \text{under~\ref{assum:maina}}, \\
2^{2(p-1)}\sup_n\norm{\vncst}^{2(p-1)}_{L^{2(p-1)}(\O)}, 	& \text{under~\ref{assum:mainb}}, \\
2^{2(p-1)}\sup_n\norm{\vncst}^{2(p-1)}_{L^{\infty}(\O)},		& \text{under~\ref{assum:mainc}},
\end{cases}
\end{align*}
and $C < +\infty$ thanks to Lemma~\ref{lem:estimatesuncstbackward}.
%
Taking $\varepsilon = \pa{\frac{C\norm{\K}_{L^{\infty,1}(\O^2)}}{(\rho_n n)^{\beta}}}^{1/2}$, we get
\begin{align*}
\P\pa{\norm{\injn\bZ}_{L^{2}(\O)} \geq \varepsilon} \leq \frac{1}{(\rho_n n)^{1-\beta}} .
\end{align*}
Plugging the latter into~\eqref{eq:unvngraphgronwall} completes the proof. 
\epf{}

\begin{rem}
Lemma~\ref{lem:intermediatebndfulldiscretebackwardgraphs} gives a deviation bound which holds with a controlled probability. On may ask if a claim with probability 1 could be afforded. A naive and straightforward approach would be to invoke the Borel-Cantelli lemma as done in~\cite[Remark~3.4(iv)]{HafieneRandom20} for the case of graphons. But this argument does not apply to the more complex setting of $L^q$-graphons given that the probability of success in the statement Lemma~\ref{lem:intermediatebndfulldiscretebackwardgraphs} does not converge sufficiently fast. This is not even possible to make faster as $\rho_n$ has to converge to $0$. Thus, it is not clear at this stage whether this is even possible to achieve or not. We leave this to a future research. 
\end{rem}
%
We finally obtain the following error bound on fully discretized problems on sparse random graphs.
\begin{theo}
Suppose that $p \in ]1, +\infty[$. Let $u$ be a solution of~\eqref{neumann} with kernel $\K$ and data $(f,g)$, and $\acc{\bu^k}_{k \in [N]}$ is the sequence generated by~\eqref{neumannfulldiscbackwardgraphs} with $\bK=\projn \K$, $\bg=\projn g$, $\bff^k = \tauh^{-1}\int_{t_{k-1}}^{t_k} \projn f(\cdot,t)dt$. Assume that $(f,g,\K)$ satisfy the assumptions of Theorem~\ref{thm:mainfulldiscretebackward}, and that those of Lemma~\ref{lem:intermediatebndfulldiscretebackwardgraphs} also hold.
\begin{enumerate}
\item For any $\beta \in ]0,1[$, with probability at least $1-(\rho_n n)^{-(1-\beta)}$,
\begin{multline}
\hspace*{-1.5cm}
\sup_{k \in [N], t \in ]t_{k-1},t_k]}\norm{\injn\bu^{k} - u(\cdot,t)}_{L^2(\O)} 
\leq \exp\pa{T/2}\Bigg(\norm{\injn\projn g - g}_{L^2(\O)} + \norm{\fncst - f}_{L^1([0,T];L^2(\O))} + CT^{1/2}(\rho_n n)^{-\beta/2} \\
+ CT^{1/2}
\begin{cases}
\tau^{1/(3-p)}+ \norm{(\K - \rho_n^{-1})_+}_{L^{\infty,2}(\O^2)} + \norm{\injn\projn \K - \K}_{L^{\infty,2}(\O^2)}	& \text{under~\ref{assum:maina}} \\
\tau^{p/(2p-1)}+ \norm{(\K - \rho_n^{-1})_+}_{L^{\infty,2}(\O^2)} + \norm{\injn\projn \K - \K}_{L^{\infty,2}(\O^2)}	& \text{under~\ref{assum:mainb}}\\
\tau^{1/(3-p)} + \norm{(\K - \rho_n^{-1})_+}_{L^{2}(\O^2)} + \norm{\injn\projn \K - \K}_{L^2(\O^2)}			& \text{under~\ref{assum:mainc} when $p \in ]1,2]$} \\
\tau + \norm{(\K - \rho_n^{-1})_+}_{L^{2}(\O^2)} +\norm{\injn\projn \K - \K}_{L^2(\O^2)}						& \text{under~\ref{assum:mainc} when $ p \geq 2 $} .
\end{cases}
\Bigg) .
\label{eq:fulldiscboundbackwardcstgraphs}
\end{multline}
for $\tau$ sufficiently small, where $C$ is positive constant that depends only on $p$, $g$, $f$ and $\K$.
\item If, moreover, $g \in L^{\infty}(\O) \cap \Lip(s,L^2(\O))$, $\K \in \Lip(s,L^2(\O^2))$, and $f \in L^1([0,T];L^{\infty}(\O)) \cap \Lip(s,L^2(\O \times [0,T]))$ then, for any $\delta \in ]0,1[$, with probability at least $1-(\rho_n n)^{-(1-\beta)}$,
\begin{multline}
\hspace*{-1cm}
\sup_{k \in [N], t \in ]t_{k-1},t_k]}\norm{\injn\bu^{k} - u(\cdot,t)}_{L^2(\O)} \leq C\exp(T/2)\Bigg((1+T^{1/2})n^{-s} + T^{1/2} \norm{(\K - \rho_n^{-1})_+}_{L^{2}(\O^2)} \\
+ T^{1/2}(\rho_n n)^{-\beta/2} + T^{1/2}\Bigg(
\begin{cases}
\tau^{\min(s,1/(3-p))}		& \text {when $p \in ]1,2]$} \\
\tau^s				& \text{ when $ p \geq 2 $} 
\end{cases}\Bigg)\Bigg).
\label{eq:fulldiscretebackwardratecstgraphs}
\end{multline}
for $\tau$ sufficiently small, where $C$ is positive constant that depends only on $p$, $g$, $f$, $\K$ and $s$, and $\norm{(\K - \rho_n^{-1})_+}_{L^{2}(\O^2)} = o(1)$. The term $\tau^s$ in the dependence on $\tau$ disappears when $f$ is time-independent. 
\end{enumerate}
\label{thm:mainfulldiscretebackwardgraphs}
\end{theo}

\bpf{}
In view of \eqref{eq:errorbnddecomp}, we shall use Theorem~\ref{thm:mainfulldiscretebackward} to bound the second term, and a bound on the first term is provided by Lemma~\ref{lem:intermediatebndfulldiscretebackwardgraphs}. Since $\injn\wedg(x,y) \leq \injn\bK(x,y) = \injn\projn\K(x,y)$, the assumptions on $\bK$ transfer to $\wedg$, and the second term of \eqref{eq:errorbnddecomp} can then be bounded using~\eqref{eq:fulldiscboundbackwardcst}, replacing $\injn\projn \K$ there by $\injn\wedg$. Observing that
\begin{align*}
\norm{\injn\wedg - \K}_{L^2(\O^2)} 
&= \norm{\min(\injn\projn\K,\rho_n^{-1}) - \K}_{L^2(\O^2)} \\
&\leq \norm{\min(\injn\projn\K,\rho_n^{-1}) - \injn\projn\K}_{L^2(\O^2)} + \norm{\injn\projn\K - \K}_{L^2(\O^2)} \\
&= \norm{(\injn\projn\K - \rho_n^{-1})_+}_{L^2(\O^2)} + \norm{\injn\projn\K - \K}_{L^2(\O^2)} \\
&\leq \norm{(\K - \rho_n^{-1})_+}_{L^2(\O^2)} + 2\norm{\injn\projn\K - \K}_{L^2(\O^2)} ,
\end{align*}
and similarly for the $L^{\infty,2}$ norm. The fact that $\norm{(\K - \rho_n^{-1})_+}_{L^{2}(\O^2)} = o(1)$ is because $\rho_n \to 0$ by the same argument as the end of the proof of Proposition~\ref{prop:sparsegraphconv}. This completes the proof.
\epf{}

\paragraph{Acknowledgement.}
This work was supported by the Normandy Region grant MoNomads.

\bibliographystyle{abbrv}
\bibliography{biblio}

\end{document}